\documentclass{amsart}

\usepackage[a4paper,hmargin=3.5cm,vmargin=4cm]{geometry}
\usepackage{amsfonts,amssymb,amscd,amstext}
\usepackage{graphicx}
\usepackage[dvips]{epsfig}
\usepackage{hyperref}

\renewcommand{\leq}{\leqslant}
\renewcommand{\geq}{\geqslant}
\newcommand{\ptl}{\partial}
\newcommand{\rr}{{\mathbb{R}}}
\newcommand{\la}{\lambda}
\newcommand{\hh}{{\mathbb{H}}}
\newcommand{\h}{\mathcal{H}}

\newcommand{\sub}{\subset}

\newcommand{\escpr}[1]{\big<#1\big>}
\newcommand{\Sg}{\Sigma} \newcommand{\sg}{\sigma}
\newcommand{\Om}{\Omega}
\newcommand{\eps}{\varepsilon}
\newcommand{\var}{\varphi}
\newcommand{\ga}{\gamma}
\newcommand{\Ga}{\Gamma}
\newcommand{\mnh}{|N_{h}|}
\newcommand{\nuh}{\nu_{h}}
\newcommand{\ric}{\text{Ric}}

\DeclareMathOperator{\divv}{div}
\DeclareMathOperator{\supp}{supp}

\newtheorem{theorem}{Theorem}[section]
\newtheorem{proposition}[theorem]{Proposition}
\newtheorem{lemma}[theorem]{Lemma}
\newtheorem{corollary}[theorem]{Corollary}

\theoremstyle{definition}

\newtheorem{remark}[theorem]{Remark}

\theoremstyle{remark}

\numberwithin{equation}{section}
\begin{document}

\bibliographystyle{amsplain}

\title[Complete stable surfaces in the Heisenberg group $\hh^1$]{the
classification of complete stable area-stationary surfaces in the
Heisenberg group $\hh^1$}

\author[A.~Hurtado]{Ana Hurtado}
\address{Departamento de Geometr\'{\i}a y Topolog\'{\i}a \\
Universidad de Granada \\ E--18071 Granada \\ Spain}
\email{ahurtado@ugr.es}

\author[M.~Ritor\'e]{Manuel Ritor\'e}
\address{Departamento de Geometr\'{\i}a y Topolog\'{\i}a \\
Universidad de Granada \\ E--18071 Granada \\ Spain}
\email{ritore@ugr.es}

\author[C.~Rosales]{C\'esar Rosales}
\address{Departamento de Geometr\'{\i}a y Topolog\'{\i}a \\
Universidad de Granada \\ E--18071 Granada \\ Spain}
\email{crosales@ugr.es}

\date{\today}

\thanks{The first author has been supported by MCyT-Feder grant
MTM2007-62344 and the Caixa-Castell\'o Foundation.  The second and
third authors have been supported by MCyT-Feder grant
MTM2007-61919 and Junta de Andaluc\'ia grant P06-FQM-01642}

\begin{abstract}
We prove that any $C^2$ complete, orientable, connected, stable
area-stationary surface in the sub-Riemannian Heisenberg group $\hh^1$
is either a Euclidean plane or congruent to the hyperbolic paraboloid
$t=xy$.
\end{abstract}

\maketitle

\thispagestyle{empty}

\section{Introduction}
\label{sec:intro}
\setcounter{equation}{0}

Minimal surfaces in Euclidean space are area-stationary, a condition
which is equivalent, by the Euler-Lagrange equation, to have mean
curvature zero.  An important question for such a variational problem
is the classification of global minimizers.  Hence is natural to
consider the second variation.  Minimal surfaces with non-negative
second variation of the area are called \emph{stable minimal
surfaces}.  It is well-known that minimal graphs are stable minimal
surfaces (in fact area-minimizing by a standard calibration argument).
A complete minimal graph must be a plane by the classical Bernstein's
Theorem \cite{bernstein}.  Bernstein result was later extended by do
Carmo and Peng \cite{dcp}, and Fischer-Colbrie and Schoen \cite{fcs},
who proved that a complete stable oriented minimal surface in $\rr^3$
must be a plane.  The proof in \cite{fcs} follows from more general
results for $3$-manifolds of non-negative scalar curvature.  Non
existence of non-orientable complete stable minimal surfaces in
$\rr^3$ has been proved by Ros \cite{ros}.

A similar analysis of the variational properties of area-minimizing
surfaces is also of great interest in some special spaces, such as the
three-dimensional Heisenberg group $\hh^1$.  This is the simplest
model of a sub-Riemannian space and of a Carnot group.  It is also the
local model of any $3$-dimensional pseudo-hermitian manifold.  For
background on $\hh^1$ we refer the reader to
Section~\ref{sec:preliminaries} and \cite{survey}.

Area-stationary surfaces of class $C^2$ in $\hh^1$ are well
understood.  It is well-known \cite{chmy}, \cite{rr2} that, outside
the singular set given by the points where the tangent plane is
horizontal, such a surface is ruled by characteristic horizontal
segments.  Moreover, based on the description of the singular set for
$t$-graphs of class $C^2$ given by Cheng, Hwang, Malchiodi and Yang
\cite{chmy}, Ritor\'e and Rosales \cite{rr2} proved that a $C^2$
surface $\Sg$ immersed in $\hh^1$ is area-stationary if and only if
its mean curvature is zero and the characteristic segments in $\Sg$
meet orthogonally the singular curves.  A similar result was
independently obtained for area-minimizing $t$-graphs by Cheng, Hwang,
and Yang \cite{chy}.  Furthermore, the classification of $C^2$
complete, connected, orientable, area-stationary surfaces with
non-empty singular set was provided in \cite{rr2}: the only examples
are, modulo congruence, non-vertical Euclidean planes, the hyperbolic
paraboloid $t=xy$, and the classical left-handed minimal helicoids.
Though some results for complete area-stationary surfaces with empty
singular set have been proved, see for example \cite[Thm.~5.4]{rr1},
\cite{chenghwang} and \cite[Prop.~6.16]{rr2}, a detailed description
of such surfaces seems far from being established.  This provides an
additional motivation for the study of second order minima of the area
in $\hh^1$.

As in the Euclidean case, we define a \emph{stable area-stationary}
surface in $\hh^1$ as a $C^2$ area-stationary surface with
non-negative second derivative of the area under compactly supported
variations.  These surfaces have been considered in
previous papers in connection with some Bernstein type problems in
$\hh^1$.  Let us describe some related works.

In \cite{chmy}, a classification of all the complete $C^2$ solutions
to the minimal surfaces equation for $t$-graphs in $\hh^1$ is given.
In \cite{rr2}, this classification was refined by showing that the
only complete area-stationary $t$-graphs are Euclidean non-vertical
planes or those congruent to the hyperbolic paraboloid $t=xy$.  By
means of a calibration argument it is also proved in \cite{rr2} that
they are all area-minimizing.

In \cite{dgn3} and \cite{bscv} the Bernstein problem for
\emph{intrinsic graphs} in $\hh^1$ was studied.  The notion of
intrinsic graph is the one used by Franchi, Serapioni and Serra
Cassano in \cite{fsscadv}.  Geometrically, an intrinsic graph is a
normal graph over some Euclidean vertical plane with respect to the
left invariant Riemannian metric $g$ in $\hh^1$ defined in
Section~\ref{sec:preliminaries}.  A $C^1$ intrinsic graph has empty
singular set.  Examples of $C^2$ complete area-stationary intrinsic
graphs different from vertical Euclidean planes were found in
\cite{dgn3}.  So a natural question is to study complete
area-minimizing intrinsic graphs.  A remarkable difference with
respect to the case of the $t$-graphs is the existence of complete
$C^2$ area-stationary intrinsic graphs which are not area-minimizing,
see \cite{dgn3}.  In \cite{bscv}, Barone, Serra Cassano and Vittone
classified complete $C^2$ area-stationary intrinsic graphs.  Then they
computed the second variation formula of the area for such graphs to
establish that the only stable ones are the Euclidean vertical planes.
An interesting calibration argument, also given in \cite{bscv}, yields
that the vertical planes are in fact area-minimizing surfaces in
$\hh^1$.

In the interesting paper \cite{dgnp}, it is proven that $C^2$ complete
stable area-stationary Euclidean graphs with empty singular set must
be vertical planes.  This is done by showing that if such a graph is
different from a vertical plane then it contains a particular example
of unstable surfaces called \emph{strict graphical strips}.  From the
geometrical point of view, a graphical strip is a $C^2$ surface given
by the union of a family of horizontal lines $L_{t}$ passing through
and filling a vertical segment so that the angle function of the
horizontal projection of $L_{t}$ is a monotonic function.  The
graphical strip is strict if the angle function is strictly monotonic.
If the angle function is constant we have a piece of a vertical plane.
We would like to remark that there are examples of complete
area-stationary surfaces with empty singular set which do not contain
a graphical strip, such as the sub-Riemannian catenoids
$t^2=\la^2\,(x^2+y^2-\la^2)$, $\la\neq 0$.  Hence the main result in
\cite{dgnp} does not apply to general surfaces.

The following natural step is to consider complete stable surfaces in
$\hh^1$.  In fact, all the aforementioned results leave open the
existence of stable examples different from intrinsic graphs or
Euclidean graphs with empty singular set.  The purpose of the present
paper is to classify complete stable area-stationary surfaces in
$\hh^1$ with empty singular set or not.  In Theorem~\ref{th:main} we
prove the following~result

\begin{quotation}
\emph{The only complete, orientable, connected, stable area-stationary
surfaces in $\hh^1$ of class $C^2$ are the Euclidean planes and the
surfaces congruent to the hyperbolic paraboloid $t=xy$.}
\end{quotation}
In particular, this result provides the classification of all the
complete $C^2$ orientable area-minimizing surfaces in $\hh^1$.

In order to prove Theorem~\ref{th:main} we compute the second derivative of the area for some compactly supported variations of a $C^2$ area-stationary surface $\Sg$ by means of Riemannian geodesics. In Theorem~\ref{th:2ndvar}, variations of a portion $\Sg'$ of the regular part of $\Sg$ in the direction of $vN+wT$, where $N$ is the unit normal to $\Sg$ and $T$ is the Reeb vector field in $\hh^1$, will be considered. Here $v$, $w$ are assumed to have compact support in $\Sg$, but not on $\Sg'$. Hence the boundary of $\Sg'$ is moving along the variation. In Proposition~\ref{prop:2ndvert} variations in the direction of $wT$ of a $C^2$ area-stationary surface $\Sg$ with singular curves of class $C^3$ will be taken. Here $w$ has compact support near the singular curves, and it is constant along the characteristic curves of $\Sg$. Both types of variations will be combined to produce global ones in Proposition~\ref{prop:2ndvarsing}.
 Second variation formulas of the area for variations \emph{supported in the
regular set} have appeared in several contexts.  In \cite{chmy}, such
a formula was obtained for $C^3$ surfaces inside a $3$-dimensional
pseudo-hermitian manifold.  In \cite{bscv}, a second variation formula
was proved for variations by intrinsic graphs of class $C^2$.  In
\cite{dgn}, it is computed the second derivative of the area associated
to a $C^2$ variation of a $C^2$ surface along Euclidean straight lines.

Once we have the second variation formula we proceed into two steps.
First we prove in Theorem~\ref{th:nosing} that a $C^2$ complete
oriented stable area-stationary surface with empty singular set must
be a vertical plane.  In fact, for such a surface $\Sg$, the second
derivative of the area for a compactly supported variation as in
Theorem~\ref{th:2ndvar} is given by
\[
\mathcal{I}(u,u)=-\int_{\Sg} u\,\mathcal{L}(u),
\]
where $u$ is the normal component of the variation, and $\mathcal{L}$
is the hypoelliptic operator on $\Sg$ given in \eqref{eq:lu}.  By
analogy with the Riemannian situation \cite{bdce} we refer to
$\mathcal{I}$ as the \emph{index form} associated to $\Sg$ and to
$\mathcal{L}$ as the \emph{stability operator} of $\Sg$.  In
Proposition~\ref{prop:stcond1} we see that the stability condition for
$\Sg$ implies that $\mathcal{I}(u,u)\geq 0$ for any $u\in C_{0}(\Sg)$
which is also $C^1$ along the characteristic lines.  Then we choose
the function $u:=|N_{h}|$, where $N$ is the Riemannian unit normal to
$\Sg$ for the left invariant Riemannian metric $g$ on $\hh^1$ defined
in Section~\ref{sec:preliminaries}, $N_{h}$ is the horizontal
projection of $N$, and the modulus is computed with respect to the
metric $g$.  We see in Proposition~\ref{lem:lnh>0} that this function
$u$ satisfies
\[
\mathcal{L}(u)\geq 0,
\]
and the inequality is strict in pieces of $\Sg$ which are not
contained inside Euclidean vertical planes.  In such a case we produce a
compactly supported non-negative function $v$ in $\Sg$ so that
inequality $\mathcal{I}(v,v)<0$ still holds.  To construct the
function $v$ we use the Jacobi vector field on $\Sg$ associated to the
family of horizontal straight lines ruling $\Sg$ and which is studied
in Lemma~\ref{lem:jacobifield}.  Observe that the function $\mnh$ is
associated to the variational vector field induced by the surfaces
equidistant to $\Sg$ in the Carnot-Carath\'eodory distance, see
\cite{arcfer}.  Hence, our construction of the test function $v$ is,
in spirit, similar to that in the Euclidean case, where the equivalent
test function is $u\equiv 1$.  Using Fischer-Colbrie's results
\cite{fc}, a stable minimal surface in $\rr^3$ is conformally a compact Riemann
surface minus a finite number of points, so that a logarithmic cut-off
function $v$ of $u\equiv 1$ has compact support and yields instability
unless the surface is a plane.  We remark that the function $\mnh$ was
already used as a test function in \cite{bscv}, \cite{dgn3} and
\cite{dgnp}.

In the second step of the proof of Theorem~\ref{th:main} we consider a
complete area-stationary surface $\Sg$ with non-empty singular set.
From the classification in \cite{rr2}, we conclude that $\Sg$ must be
a non-vertical plane, congruent to the hyperbolic paraboloid $t=xy$, or
congruent to a left-handed helicoid, see Proposition~\ref{prop:basic2}
for a precise statement.  The first two types of surfaces are
$t$-graphs and then they are area-minimizing by a calibration argument
\cite{rr2}. For the third type we will combine our second variation formulas in Theorem~\ref{th:2ndvar} and Proposition~\ref{prop:2ndvert} to produce the stability inequality $\mathcal{Q}(u)\geq 0$, where $\mathcal{Q}$ is the quadratic form defined in \eqref{eq:Qu}. The construction of appropriate test functions with $\mathcal{Q}(u)<0$ will prove the instability of the helicoids. It is interesting to observe that $\mathcal{Q}(u)\geq 0$ for functions $u$ with support in the regular part of the helicoids.

In the Heisenberg groups $\hh^n$, with $n\geq 5$, there is no
counterpart to Theorem~\ref{th:nosing}, as some examples have been
constructed in \cite{bscv} of complete area-minimizing intrinsic
graphs different from Euclidean vertical hyperplanes.  For $n=2,3,4$
it is still unknown if similar examples can be obtained.

We would like to mention that examples of area-minimizing surfaces in
$\hh^1$ with low Euclidean regularity have been obtained in
\cite{chy}, \cite{pauls-regularity}, \cite{r2} and \cite{mscv}.  Hence
our results are optimal in the class of $C^2$ area-stationary
surfaces.

Finally, the techniques in this paper can be employed to prove
classification results for complete stable area-stationary surfaces
under a volume constraint in the first Heisenberg group \cite{rcmc},
and inside the sub-Riemannian three-sphere \cite{hr2}.

We have organized this paper as follows: the next section contains
some background material in several subsections.  In the third one we
recall known facts about area-stationary surfaces and we compute
second variation formulas for the area.  The fourth and fifth
sections treat complete stable surfaces without and with singular
points, respectively.  In the sixth section we state and prove the
main result.

After the distribution of this paper we were informed by Prof. Nicola Garofalo that Theorem~\ref{th:nosing} was proven, for the case of embedded surfaces, by Danielli, Garofalo, Nhieu and Pauls in late 2006, \cite{dgnp-stable}.

\section{Preliminaries}
\label{sec:preliminaries}
\setcounter{equation}{0}

In this section we gather some previous results that will be used
throughout the paper.  We have organized it in several parts.

\subsection{The Heisenberg group}
\label{subsec:hg}
The \emph{Heisenberg group} $\hh^1$ is the Lie group $(\rr^3,*)$,
where the product $*$ is defined, for any pair of points $[z,t]$,
$[z',t']\in\rr^3\equiv\mathbb{C}\times\rr$, by
\[
[z,t]*[z',t']:=[z+z',t+t'+\text{Im}(z\overline{z}')], \qquad (z=x+iy).
\]
For $p\in\hh^1$, the \emph{left translation} by $p$ is the
diffeomorphism $L_p(q)=p*q$.  A basis of left invariant vector fields
(i.e., invariant by any left translation) is given by
\begin{equation*}
X:=\frac{\ptl}{\ptl x}+y\,\frac{\ptl}{\ptl t}, \qquad
Y:=\frac{\ptl}{\ptl y}-x\,\frac{\ptl}{\ptl t}, \qquad
T:=\frac{\ptl}{\ptl t}.
\end{equation*}
The \emph{horizontal distribution} $\mathcal{H}$ in $\hh^1$ is the
smooth planar distribution generated by $X$ and $Y$.  The
\emph{horizontal projection} of a tangent vector $U$ onto
$\mathcal{H}$ will be denoted by $U_{h}$.  A vector field $U$ is
\emph{horizontal} if $U=U_h$.

We denote by $[U,V]$ the Lie bracket of two $C^1$ vector fields $U$
and $V$ on $\hh^1$.  Note that $[X,T]=[Y,T]=0$, while $[X,Y]=-2T$, so
that $\mathcal{H}$ is a bracket-generating distribution.  Moreover, by
Frobenius theorem we have that $\mathcal{H}$ is nonintegrable.  The
vector fields $X$ and $Y$ generate the kernel of the (contact)
$1$-form $\omega:=-y\,dx+x\,dy+dt$.

\subsection{The left invariant metric}
\label{subsec:g}
We shall consider on $\hh^1$ the Riemannian metric
$g=\escpr{\cdot\,,\cdot}$ so that $\{X,Y,T\}$ is an orthonormal basis
at every point.  The restriction of $g$ to $\h$ coincides with the
usual sub-Riemannian metric in $\hh^1$.  Let $D$ be the Levi-Civita
connection associated to $g$.  From Koszul formula and the Lie bracket
relations we get
\begin{alignat}{2}
\notag
D_{X}X&=0, \qquad \ \ \ \, D_{Y}Y=0, \qquad \,D_{T}T=0, \\
\label{eq:christoffel}
D_{X}Y&=-T, \qquad   \, D_{X}T=Y, \qquad   D_{Y}T=-X, \\
\notag D_{Y}X&=T, \qquad \ \ \,D_{T}X=Y, \qquad D_{T}Y=-X.
\end{alignat}
For any tangent vector $U$ on $\hh^1$ we define $J(U):=D_UT$.  Then we
have $J(X)=Y$, $J(Y)=-X$ and $J(T)=0$, so that $J^2=-\text{Id}$
when restricted to $\h$.  It is also clear~that
\begin{equation}
\label{eq:conmute}
\escpr{J(U),V}+\escpr{U,J(V)}=0,
\end{equation}
for any pair of tangent vectors $U$ and $V$.  The involution
$J:\h\to\h$ together with the $1$-form $\omega=-y\,dx+x\,dy+dt$,
provides a pseudo-hermitian structure on $\hh^1$, see
\cite[Sect.~6.4]{blair}.

Let $R$ be the Riemannian curvature tensor of $g$ defined for tangent
vectors $U,V,W$ by
\[
R(U,V)W=D_{V}D_{U}W-D_{U}D_{V}W+D_{[U,V]}W.
\]
From \eqref{eq:christoffel} and the Lie bracket relations we can
obtain the following identities
\begin{align}
\notag
R(X,Y)\,X&=-3Y, & R(X,Y)\,Y&=3X, & R(X,Y)\,T&=0,
\\
\label{eq:curvature}
R(X,T)\,X&=T, & R(X,T)\,Y&=0, & R(X,T)\,T&=-X,
\\
\nonumber
R(Y,T)\,X&=0, & R(Y,T)\,Y&=T, & R(Y,T)\,T&=-Y.
\end{align}
We denote by $\text{Ric}$ the Ricci curvature in $(\hh^1,g)$ defined,
for any pair of tangent vectors $U$ and $V$, as the trace of the map
$W\mapsto R(U,W)V$.  These equalities can be checked by taking into
account \eqref{eq:curvature}
\begin{align}
\label{eq:ricci}
\ric(X,Y)&=0, & \ric(X,T)&=0, & \ric(Y,T)&=0,
\\
\notag
\ric(X,X)&=-2, & \ric(Y,Y)&=-2, & \ric(T,T)&=2.
\end{align}

\subsection{Horizontal curves and Carnot-Carath\'eodory distance}
Let $\ga:I\to\hh^1$ be a piecewise $C^1$ curve defined on a compact
interval $I\sub\rr$.  The \emph{length} of $\ga$ is the usual
Riemannian length $L(\ga):=\int_{I}|\dot{\ga}(\eps)|\,d\eps$, where
$\dot{\ga}$ is the tangent vector of $\ga$.  A \emph{horizontal curve}
$\ga$ in $\hh^1$ is a $C^1$ curve whose tangent vector always lies in
the horizontal distribution.  For two given points in $\hh^1$ we can
find, by Chow's connectivity theorem \cite[Sect.~1.2.B]{gromov-cc}, a
horizontal curve joining these points.  The
\emph{Carnot-Carath\'eodory distance} $d_{cc}$ between two points in
$\hh^1$ is defined as the infimum of the length of horizontal curves
joining the given points.  The topology associated to $d_{cc}$
coincides with the usual topology in $\rr^3$, see
\cite[Cor.~2.6]{andre}.

\subsection{Geodesics and Jacobi fields in $(\hh^1,g)$}
\label{subsec:geo}
A \emph{geodesic} in $(\hh^1,g)$ is a $C^2$ curve $\ga$ such that the
covariant derivative of the tangent vector field $\dot{\ga}$ vanishes
along $\ga$.

Let $\ga(s)=(x(s),y(s),t(s))$.  Dots will indicate derivatives with
respect to $s$.  We write
$\dot{\ga}=\dot{x}\,X+\dot{y}\,Y+(\dot{t}-\dot{x}y+x\dot{y})\,T$.
Then $\ga$ is a geodesic in $(\hh^1,g)$ if and only if
\begin{align*}
\ddot{x}&=2\,\escpr{\dot{\ga},T}\,\dot{y}, \\
\ddot{y}&=-2\,\escpr{\dot{\ga},T}\,\dot{x}, \\
\frac{d}{ds}&\,\escpr{\dot{\ga},T}=0.
\end{align*}
Let $\la$ be the constant
$\dot{t}-\dot{x}y+x\dot{y}=\escpr{\dot{\ga},T}$.  An easy integration
shows that the geodesic with initial conditions
$(x(0),y(0),t(0))=(x_{0},y_{0},t_{0})$ and
$(\dot{x}(0),\dot{y}(0),\dot{t}(0))=(A,B,C)$ is given by
\begin{align}
\nonumber
x(s)&=x_{0}+As\,f(2\la s)+Bs\,g(2\la s),
\\
\label{eq:eqgeo}
y(s)&=y_{0}-As\,g(2\la s)+Bs\, f(2\la s),
\\
\nonumber
t(s)&=t_{0}+\la s+(A^2+B^2) s^2\,h(2\la s)+(Ax_{0}+By_{0})s\,g(2\la s)
\\
\nonumber
&\hspace{1em}+(Ay_{0}-Bx_{0})s\,f(2\la s),
\end{align}
where $f$, $g$ and $h$ are the real analytic functions
\[
f(x):=
\begin{cases}
\displaystyle\frac{\sin(x)}{x}, \! &x\neq 0
\\
1, \! &x=0
\end{cases},
\quad
g(x):=
\begin{cases}
\displaystyle\frac{1-\cos(x)}{x}, \! &x\neq 0
\\
0, \, &x=0
\end{cases},
\quad
h(x):=
\begin{cases}
\displaystyle\frac{x-\sin(x)}{x^2}, \! &x\neq 0
\\
0, \, &x=0
\end{cases}.
\]
In particular, we have
\begin{equation}
\label{eq:horgeo}
\exp_{p}(sv)=p+sv,\qquad\text{for } \, p\in\hh^1 \ \text{and }
\, v\in\mathcal{H}_{p} \ \text{or }\,v\,||\,T_{p},
\end{equation}
which is a horizontal or vertical straight line. Here $\exp_{p}$ denotes the
exponential map of $(\hh^1,g)$ at~$p$.

In the next result we construct Riemannian Jacobi fields associated to
$C^1$ families of Riemannian geodesics.

\begin{lemma}
\label{lem:jacobi}
Let $\alpha:I\to\hh^1$ be a $C^1$ curve defined on some open interval
$I\subseteq\rr$.  For any $C^1$ vector field $U$ along $\alpha$ we
consider the map $F:I\times\rr\to\hh^1$ given by
$F(\eps,s):=\exp_{\alpha(\eps)}(s\,U_{\alpha(\eps)})$.  Then, the
variational vector field $V_{\eps}(s):=(\ptl F/\ptl\eps)(\eps,s)$ is
$C^\infty$ along the geodesic $\ga_{\eps}(s):=F(\eps,s)$.  As a
consequence, $[\dot{\ga}_{\eps},V_{\eps}]=0$ and $V_{\eps}$ satisfies
the Jacobi equation
\begin{equation}
\label{eq:jaceq1}
V_{\eps}''+R(\dot{\ga}_{\eps},V_{\eps})\dot{\ga}_{\eps}=0,
\end{equation}
where the prime $'$ denotes the covariant derivative along the geodesic
$\ga_{\eps}$.  Moreover, if $\ga_{\eps}$ is a horizontal straight
line, then
\begin{equation}
\label{eq:jaceq2}
V_{\eps}''-
3\escpr{V_{\eps},J(\dot{\ga}_{\eps})}\,J(\dot{\ga}_{\eps})+
|\dot{\ga}_{\eps}|^2\,\escpr{V_{\eps},T}\,T=0.
\end{equation}
\end{lemma}

\begin{remark}
\label{re:riem}
The classical proofs in Riemannian geometry of
$[\dot{\ga}_{\eps},V_{\eps}]=0$ and the fact that $V_{\eps}$ satisfies
the Jacobi equation do not apply directly in our setting since we only
suppose that $F$ is a $C^1$ map.
\end{remark}

\begin{proof}[Proof of Lemma~\ref{lem:jacobi}]
Let $(x_{0}(\eps),y_{0}(\eps),t_{0}(\eps))$ and
$(A(\eps),B(\eps),C(\eps))$ be the Euclidean coordinates of
$\alpha(\eps)$ and $U_{\alpha(\eps)}$, respectively.  By using the
expression of the Riemannian geodesics in \eqref{eq:eqgeo}, we see
that the map $F(\eps,s)$ can be written as
\begin{align*}
x(\eps,s)&=x_{0}(\eps)+A(\eps)s\,f(2\la(\eps)s)+B(\eps)s\,g(2\la(\eps) s),
\\
y(\eps,s)&=y_{0}(\eps)-A(\eps)s\,g(2\la(\eps) s)+B(\eps)s\, f(2\la(\eps) s),
\\
t(\eps,s)&=t_{0}(\eps)+\la(\eps) s+(A^2+B^2)(\eps) s^2\,h(2\la(\eps)
s)+(A(\eps)x_{0}(\eps)+B(\eps)y_{0}(\eps))s\,g(2\la(\eps) s)
\\
&\hspace{1em}+(A(\eps)y_{0}(\eps)-B(\eps)x_{0}(\eps))s\,f(2\la(\eps) s),
\end{align*}
where $\la(\eps):=C(\eps)-A(\eps)y_{0}(\eps)+B(\eps)x_{0}(\eps)$.
Observe that the functions $x_{0}(\eps)$, $y_{0}(\eps)$,
$t_{0}(\eps)$, $A(\eps)$, $B(\eps)$, $C(\eps)$ and $\la(\eps)$ are
$C^1$.  A direct computation of $(\ptl F/\ptl\eps)(\eps,s)$ shows that
$V_{\eps}(s)$ is $C^\infty$ along the geodesic $\ga_{\eps}(s)$.

On the other hand, we can check that for all $k\in\mathbb{N}$ and any
of the Euclidean components $\phi(\eps,s)$ of $F(\eps,s)$, the partial
derivatives $\ptl^{k+1}\phi/\ptl\eps\,\ptl^k s$ exist and are continuous
functions.  In particular, it follows from the classical Schwarz's
theorem that $\ptl^2\phi/\ptl\eps\ptl s=\ptl^2\phi/\ptl s\ptl\eps$ and
$\ptl^3\phi/\ptl\eps\ptl s^2=\ptl^3\phi/\ptl s\ptl\eps\ptl s$.  Now,
the classical proofs in \cite[p.~68 and p.~111]{dcriem} can be
traced to prove that $[\dot{\ga}_{\eps},V_{\eps}]=0$ and that
$V_{\eps}$ satisfies the Jacobi equation.  Finally, to get
\eqref{eq:jaceq2} from \eqref{eq:jaceq1} it suffices to use
\eqref{eq:curvature} to obtain
$R(w,v)w=-3\escpr{v,J(w)}J(w)+|w|^2\escpr{v,T}T$ provided $w$ is a
horizontal vector.
\end{proof}

\subsection{Geometry of surfaces in $\hh^1$}
\label{subsec:surfaces}
Unless explicitly stated we shall consider surfaces with empty
boundary.
Let $\Sg$ be a $C^1$ surface immersed in $\hh^1$.  The \emph{singular
set} $\Sg_0$ consists of those points $p\in\Sg$ for which the tangent
plane $T_p\Sg$ coincides with $\h_{p}$.  As $\Sg_0$ is closed and has
empty interior in $\Sg$, the \emph{regular set} $\Sg-\Sg_0$ of $\Sg$
is open and dense in $\Sg$.  It was proved in \cite[Lem.~1]{d2}, see
also \cite[Thm.~1.2]{balogh}, that, for a $C^2$ surface, the Hausdorff
dimension of $\Sg_{0}$ with respect to the Riemannian distance on
$\hh^1$ is less than or equal to one.  In particular, the Riemannian
area of $\Sg_{0}$ vanishes.  If $N$ is a unit normal vector to $\Sg$
in $(\hh^1,g)$, then we can describe the singular set as
$\Sg_{0}=\{p\in\Sg;N_h(p)=0\}$, where $N_{h}=N-\escpr{N,T}T$.  In the
regular part $\Sg-\Sg_0$, we can define the \emph{horizontal Gauss
map} $\nu_h$ and the \emph{characteristic vector field} $Z$, by
\begin{equation}
\label{eq:nuh}
\nu_h:=\frac{N_h}{|N_h|}, \qquad Z=J(\nuh).
\end{equation}
As $Z$ is horizontal and orthogonal to $\nu_h$, we conclude that $Z$
is tangent to $\Sg$.  Hence $Z_{p}$ generates $T_{p}\Sg\cap\h_{p}$.
The integral curves of $Z$ in $\Sg-\Sg_0$ will be called
$(\!\emph{oriented})$ \emph{characteristic curves} of $\Sg$.  They are
both tangent to $\Sg$ and horizontal.  If we define
\begin{equation}
\label{eq:ese}
S:=\escpr{N,T}\,\nu_h-|N_h|\,T,
\end{equation}
then $\{Z_{p},S_{p}\}$ is an orthonormal basis of $T_p\Sg$ whenever
$p\in\Sg-\Sg_0$.  Moreover, for any $p\in\Sg-\Sg_{0}$ we have the
orthonormal basis of $T_{p}\hh^1$ given by $\{Z_{p},(\nuh)_{p},T_{p}
\}$.  From here we deduce the following identities on $\Sg-\Sg_{0}$
\begin{equation}
\label{eq:relations}
|N_{h}|^2+\escpr{N,T}^2=1, \qquad (\nu_{h})^\top=\escpr{N,T}\,S, \qquad
T^\top=-|N_{h}|\,S,
\end{equation}
where $U^\top$ stands for the projection of a vector field $U$ onto
the tangent plane to $\Sg$.

Given a $C^1$ immersed surface $\Sg$ with a unit normal vector $N$, we
define the \emph{area} of $\Sg$ by
\begin{equation}
\label{eq:area}
A(\Sg):=\int_{\Sg}|N_{h}|\,d\Sg,
\end{equation}
where $d\Sg$ is the Riemannian area element on $\Sg$.  If $\Sg$ is a
$C^2$ surface bounding a set $\Om$, then $A(\Sg)$ coincides with
all the notions of perimeter of $\Om$ and area of $\Sg$ introduced by
other authors, see \cite[Prop.~2.14]{fssc}, \cite[Thm.~5.1]{msc}
and \cite[Cor.~7.7]{fssc}.

Finally, for a $C^2$ immersed surface $\Sg$ with a unit normal vector
$N$, we denote by $B$ the Riemannian shape operator of $\Sg$ with
respect to $N$.  It is defined for any vector $W$ tangent to $\Sg$ by
$B(W)=-D_{W}N$.  The Riemannian mean curvature of $\Sg$ is
$-2H_{R}=\divv_{\Sg}N$, where $\divv_{\Sg}$ denotes the Riemannian
divergence relative to $\Sg$.

\subsection{Isometries and dilations}
\label{subsec:iso}
By a \emph{horizontal isometry} of $\hh^1$ we mean an isometry of
$(\hh^1,g)$ leaving invariant the horizontal distribution.  These
isometries preserve the area defined in \eqref{eq:area}.  Examples of
such isometries are the left translations and the Euclidean rotations
about the $t$-axis.  We say that two surfaces $\Sg_{1}$ and $\Sg_{2}$
are \emph{congruent} if there is a horizontal isometry $\phi$ such
that $\phi(\Sg_{1})=\Sg_{2}$.

In the Heisenberg group $\hh^1$ there is a one-parameter group of
$C^\infty$ \emph{dilations} $\{\delta_\la\}_{\la\in\rr}$ given in
coordinates $(x,y,t)$ by
\begin{equation}
\label{eq:dilations}
\delta_{\la}(x,y,t)=(e^\la x,e^\la y,e^{2\la}t).
\end{equation}
From \eqref{eq:dilations} it is easy to check that any $\delta_{\la}$
preserves the horizontal and the vertical distributions.  The
behaviour of the area with respect to $\delta_{\la}$ is contained in
the formula
\begin{equation}
\label{eq:areadil}
A(\delta_{\la}(\Sg))=e^{3\la}A(\Sg).
\end{equation}
For a proof of \eqref{eq:areadil} see \cite[Proof of Thm.~4.12]{rr2}.

\subsection{A weak Riemannian divergence theorem}
\label{subsec:div}
Let $\Sg$ be a $C^2$ Riemannian surface.  For any
integer $r\geq 1$ we denote by $C^r_{0}(\Sg)$ and $C^r(\Sg)$ the
spaces of functions of class $C^r$ with or without compact support in
$\Sg$.  For $r\geq 1$ let $L^r(\Sg)$ be the corresponding space of
integrable functions with respect to the Riemannian measure $d\Sg$.
Let $U$ be a $C^1$ tangent vector field on $\Sg$.  Given a continuous
function $f$ on $\Sg$, a continuous vector field $V$ on $\Sg$, and a
point $p\in \Sg$, we define $U_{p}(f)=(f\circ\alpha)'(0)$ and
$(D_{U}V)(p)=V_{\alpha(s)}'(0)$.  Here $\alpha$ is the integral curve
of $U$ with $\alpha(0)=p$, while the primes denote derivatives of
functions depending on $s$ and covariant derivatives along
$\alpha(s)$.  We say that $f$ and $V$ are $C^1$ in the $U$-direction
if $U(f)$ and $D_{U}V$ are well defined and they are continuous on
$\Sg$.  We also set
\begin{equation}
\label{eq:divv}
\divv_{\Sg}(fU):=f\divv_{\Sg}U+U(f),
\end{equation}
where $\divv_{\Sg}U$ stands for the Riemannian divergence of $U$.  Note
that these definitions coincide with the classical ones when $f\in
C^1(\Sg)$ and $V$ is a $C^1$ vector field on $\Sg$.  In the same way we
can introduce derivatives of higher order in the $U$-direction.

Now we extend the classical Riemannian divergence theorem in $\Sg$ to
certain vector fields with compact support which are not $C^1$ on $\Sg$.
First we need an approximation result.

\begin{lemma}
\label{lem:approx}
Let $\Sg$ be a $C^2$ Riemannian surface.  Consider a
$C^1$ tangent vector field $U$ on $\Sg$ such that $U_{p}\neq 0$ for any
$p\in \Sg$.  Then, for any function $f\in C_{0}(\Sg)$ which is also $C^1$
in the $U$-direction, there is a compact set $K\subseteq \Sg$ and a
sequence of functions $\{f_{\eps}\}_{\eps>0}$ in $C^1_{0}(\Sg)$ such
that the supports of $f$ and $f_{\eps}$ are contained in $K$ for any
$\eps>0$, and
\begin{itemize}
\item[(i)] $\{f_{\eps}\}\to f$ in $L^r(\Sg)$ for any integer $r\geq 1$,
\item[(ii)] $\{U(f_{\eps})\}\to U(f)$ in $L^r(\Sg)$ for any integer $r\geq 1$.
\end{itemize}
\end{lemma}

\begin{proof}
Let $p\in \Sg$.  By using the local flow of $U$ in $\Sg$ and that
$U_{p}\neq 0$, we can find a local $C^1$ chart $(D,\phi=(x,y))$ of
$\Sg$ around $p$ such that $K=\overline{D}$ is compact and the
restriction of $U$ to $D$ coincides with the basic vector field
$\ptl_{y}$.  This means that $U(h)=(\ptl(h\circ\phi^{-1})/\ptl
y)\circ\phi$ for any function $h$ which is $C^1$ in the $U$-direction.
To finish the proof it suffices, by a standard partition of unity
argument, to prove the claim when the support of $f$ is contained in
$D$.  Let $D'=\phi(D)$ and $g=f\circ\phi^{-1}$.  We have $g\in
C_{0}(D')$ and $\ptl g /\ptl y=U(f)\circ\phi^{-1}\in C_{0}(D')$.  From
the standard regularization by convolution in $\rr^2$, see for
instance \cite[Sect.~4.2.1]{evans}, we can find a sequence
$\{g_{\eps}\}_{\eps>0}$ in $C^\infty_{0}(\rr^2)$ such that
$\{g_{\eps}\}\to g$ and $\{\ptl g_{\eps}/\ptl y\}\to \ptl g/\ptl y$
uniformly in $\rr^2$, while the supports of $g_{\eps}$ are contained
in $D'$ for any $\eps>0$.  It follows that the family
$\{f_{\eps}\}_{\eps>0}$ with $f_{\eps}=g_{\eps}\circ\phi$ satisfies
$\{f_{\eps}\}\to f$ and $\{U(f_{\eps})\}\to U(f)$ uniformly in $D$,
while the support of $f_{\eps}$ is contained in $D\sub K$ for any
$\eps>0$.  Clearly $\{f_{\eps}\}_{\eps>0}$ proves the lemma.
\end{proof}

\begin{lemma}
\label{lem:div}
Let $\Sg$ be a $C^2$ Riemannian surface.  Consider a
$C^1$ tangent vector field $U$ on $\Sg$ such that $U_{p}\neq 0$ for any
$p\in \Sg$.  Then, for any $f\in C_{0}(\Sg)$ which is also $C^1$ in the
$U$-direction, we have
\[
\int_{\Sg}\divv_{\Sg}(fU)\,d\Sg=0.
\]
\end{lemma}

\begin{proof}
By definition \eqref{eq:divv} it follows that $\divv_{\Sg}(fU)\in
L^1(\Sg)$ since $f$ has compact support and $U(f)$ is continuous.  By
Lemma~\ref{lem:approx} we can find a sequence $\{f_{\eps}\}_{\eps>0}$
in $C^1_{0}(\Sg)$ such that $\{f_{\eps}\}\to f$ and
$\{U(f_{\eps})\}\to U(f)$ in $L^1(\Sg)$, while the supports of
$f_{\eps}$ and $f$ are contained in the same compact set $K\subseteq
\Sg$ for any $\eps>0$.  In particular, we deduce
$\{f_{\eps}\divv_{\Sg}U\}\to f\divv_{\Sg}U$ in $L^1(\Sg)$ since
$\divv_{\Sg}U$ is continuous.  By using the Riemannian divergence
theorem for $C^1$ vector fields with compact support, we obtain
\[
0=\int_{\Sg}\divv_{\Sg}(f_{\eps}U)\,d\Sg=\int_{\Sg}f_{\eps}
\divv_{\Sg}U\,d\Sg+\int_{\Sg}U(f_{\eps})\,d\Sg,\qquad\eps>0.
\]
Letting $\eps\to 0$ in the previous equality the claim is proven.
\end{proof}

\section{Stable surfaces. Second variation formulas of the area}
\label{sec:2ndvar}

In this section we define stable surfaces and we show that they
satisfy an analytical inequality by means of a second variation
formula for the area functional defined in \eqref{eq:area}.  We first
introduce the appropriate variational background.

Let $\Sg$ be a $C^2$ oriented surface immersed in $\hh^1$ with
singular set $\Sg_{0}$.  By a \emph{variation} of $\Sg$ we mean a
$C^1$ map $\varphi:I\times\Sg\to\hh^1$, where $I$ is an open interval
containing the origin, satisfying the following properties:
\begin{itemize}
\item[(i)] $\var(0,p)=p$ for any $p\in\Sg$,
\item[(ii)] The set $\Sg_{s}=\{\varphi(s,p);\,p\in\Sg\}$ is a $C^1$
surface immersed in $\hh^1$ for any $s\in I$,
\item[(iii)] The map $\varphi_{s}:\Sg\to\Sg_{s}$ given
by $\varphi_{s}(p)=\varphi(s,p)$ is a diffeomorphism for any $s\in I$.
\end{itemize}
We say that the variation is \emph{compactly supported} if there is a
compact set $K\subseteq\Sg$ such that $\varphi_{s}(p)=p$ for any $s\in
I$ and $p\in\Sg-K$.  If, in addition, the set $K$ is contained inside
$\Sg-\Sg_{0}$ then the variation is \emph{nonsingular}.  The area
functional associated to the variation is $A(s):=A(\Sg_{s})$.  Note
that only the deformation over the compact set $K$ contributes to the
change of area.  We say that $\Sg$ is \emph{area-stationary} if
$A'(0)=0$ for any compactly supported variation.  We say that $\Sg$ is
\emph{stable} (resp.  \emph{stable under non-singular variations}) if
it is area-stationary and $A''(0)\geq 0$ for any compactly supported
(resp.  non-singular) variation of $\Sg$.  Finally by an
\emph{area-minimizing} surface in $\hh^1$ we mean a $C^2$ orientable surface
$\Sg$ such that any compact region $M\sub\Sg$ satisfies $A(M)\leq
A(M')$ for any other $C^1$ compact surface $M'$ in $\hh^1$ with $\ptl
M=\ptl M'$.  Clearly any area-minimizing surface is stable.

\begin{remark}
Consider a $C^1$ vector field $U$ with compact support on $\Sg$.  For
any $s\in\rr$ we denote $\varphi_{s}(p)=\exp_{p}(s U_{p})$, where
$\exp_{p}$ is the exponential map of $(\hh^1,g)$ at $p$.  It is easy
to see that, for $s$ small enough, $\{\var_{s}\}_{s}$ defines a
compactly supported variation of $\Sg$.  In case the support of $U$ is
contained in $\Sg-\Sg_{0}$ then the induced variation is nonsingular.
This was the point of view used in \cite{rr2} to define variations of
a $C^2$ surface.  In particular, our notion of area-stationary surface
implies the one introduced in \cite[Sect.~4]{rr2}.
\end{remark}

It is clear that stability is preserved under left translations and
vertical rotations since they are horizontal isometries in $\hh^1$.
In the next result we prove that any dilation $\delta_{\la}$ as
defined in \eqref{eq:dilations} satisfies the same property.

\begin{lemma}
\label{lem:dilinv}
Let $\Sg$ be a $C^2$ immersed oriented surface in $\hh^1$.  Then $\Sg$
is stable $($resp.  stable under non-singular variations$)$ if and
only if the same holds for $\delta_{\la}(\Sg)$.
\end{lemma}

\begin{proof}
Let $\Sg_{\la}=\delta_{\la}(\Sg)$.  Take a compactly supported
variation $\{\var_{s}\}_{s\in I}$ of $\Sg_{\la}$.  By using that the
family of dilations is a one-parameter group of diffeomorphisms we can
see that $\{\psi_{s}\}_{s\in I}$ with
$\psi_{s}=\delta_{-\la}\circ\var_{s}\circ\delta_{\la}$ provides a
compactly supported variation of $\Sg$.  Moreover, the variation
$\{\psi_{s}\}_{s\in I}$ is nonsingular if and only if
$\{\var_{s}\}_{s\in I}$ is nonsingular.  By \eqref{eq:areadil} we get
\[
A(\Sg_{s})=A(\psi_{s}(\Sg))=A(\delta_{-\la}((\Sg_{\la})_{s}))
=e^{-3\la} A((\Sg_{\la})_{s}).
\]
From here it is easy to deduce that if $\Sg$ is stable (resp.  stable
under non-singular variations) then the same holds for $\Sg_{\la}$.
To prove the reverse statement it suffices to change the roles of
$\Sg$ and $\Sg_{\la}$.
\end{proof}

\subsection{Area-stationary surfaces}
In this part of the section we gather some facts about area-stationary
surfaces in $\hh^1$ that will be useful in the sequel.

Let $\Sg$ be a $C^2$ immersed surface in $\hh^1$ with a unit normal
vector $N$.  We define the \emph{mean curvature} of $\Sg$ as in
\cite{rr1} and \cite{rr2}, by the equality
\begin{equation}
\label{eq:mc}
-2H(p)=(\divv_{\Sg}\nu_{h})(p),\qquad p\in\Sg-\Sg_{0},
\end{equation}
where $\nuh$ is the horizontal Gauss map defined in \eqref{eq:nuh} and
$\divv_{\Sg}U$ stands for the divergence relative to $\Sg$ of a $C^1$
vector field $U$.  We say that $\Sg$ is a \emph{minimal surface} if
the mean curvature vanishes on $\Sg-\Sg_{0}$.

In the following proposition we recall some features about
area-stationary and minimal surfaces in $\hh^1$ involving the
structure of the regular and the singular set, see
\cite[Sect.~3]{chmy}, \cite[Sect.~4]{rr2} and the references therein.
Similar results also hold in other sub-Riemannian spaces, see \cite{hp1} and
\cite{hr1}.

\begin{proposition}
\label{prop:varprop}
Let $\Sg$ be a $C^2$ immersed oriented minimal surface in $\hh^1$ with
singular set $\Sg_{0}$. Then we have
\begin{itemize}
\item[(i)] Any characteristic curve of $\Sg$ is a segment of a
horizontal straight line.  \item[(ii)] $\Sg_{0}$ consists of isolated
points and $C^1$ curves with non-vanishing tangent vector $($singular
curves$)$.  \item[(iii)] If $\Ga$ is a singular curve and $p\in\Ga$,
then there is a neighborhood $B$ of $p$ in $\Sg$ such that $B-\Ga$ is
the union of two disjoint domains $B^+$ and $B^-$ contained in
$\Sg-\Sg_{0}$.  Moreover, the vector fields $Z$ and $\nuh$ extend
continuously to $p$ from $B^+$ and $B^-$ in such a way that
$Z_{p}^+=-Z_{p}^-$ and $(\nuh)^+_p=-(\nuh)^-_p$.
\item[(iv)] If $\Sg$ is any $C^2$ immersed
oriented surface, then $\Sg$ is area-stationary if and only if $\Sg$ is
minimal and the characteristic curves meet orthogonally the
singular curves.
\end{itemize}
\end{proposition}

Now we prove a regularity result for minimal surfaces in $\hh^1$.
Given a $C^2$ surface $\Sg$ in $\hh^1$ with unit normal vector $N$, it
is clear that the vector field $D_{Z}N$ is well defined on
$\Sg-\Sg_{0}$ and it is continuous.  By using the ruling property of
minimal surfaces in Proposition~\ref{prop:varprop} (i) we can obtain
more regularity for $N$ in the $Z$-direction.

\begin{lemma}
\label{lem:reg}
Let $\Sg$ be a $C^2$ immersed oriented surface in $\hh^1$.  If $\Sg$
is minimal then, in $\Sg-\Sg_{0}$, the normal vector $N$ is $C^\infty$
in the direction of the characteristic field $Z$.
\end{lemma}

\begin{proof}
Take $p\in\Sg-\Sg_{0}$.  Let $\ga$ be the characteristic curve through
$p$.  Consider a $C^1$ curve
$\alpha:(-\eps_{0},\eps_{0})\to\Sg-\Sg_{0}$ transverse to $\ga$ with
$\alpha(0)=p$.  Define $F(\eps,s):=\alpha(\eps)+s\,Z_{\alpha(\eps)}$.
By using \eqref{eq:horgeo} and Lemma~\ref{lem:jacobi} we get that
$V(s):=(\ptl F/\ptl\eps)(0,s)$ is a $C^\infty$ Jacobi field along
$\ga$.  Since both $\dot{\ga}(s)$ and $V(s)$ are $C^\infty$ and
linearly independent for $s$ small enough, the unit normal $N$ to
$\Sg$ along $\ga$ is given by
\[
N=\pm\,\frac{\dot{\ga}\times V}
{\big|\dot{\ga}\times V\big|},
\]
where $\times$ denotes the cross product in $(\hh^1,g)$.  We conclude
that $N$ is a $C^\infty$ vector field along $\ga$.
\end{proof}

\subsection{Second variation of the area}
\label{subsec:2ndvar}
In this part of the section we provide some formulas for the second
derivative of the area functional associated to some variations of an
area-stationary surface.  We first give some preliminary computations.

\begin{lemma}
Let $\Sg\subset\hh^1$ be a $C^2$ immersed surface with unit normal
vector $N$ and singular set $\Sg_{0}$.  Consider a point
$p\in\Sg-\Sg_0$, the horizontal Gauss map $\nu_h$ and the
characteristic field $Z$ defined in \eqref{eq:nuh}.  For any $v\in
T_{p}\hh^1$ we have
\begin{align}
\label{eq:dvnh}
D_{v}N_{h}&=(D_{v}N)_{h}-\escpr{N,T}\,J(v)
-\escpr{N,J(v)}\,T,
\\
\label{eq:vmnh}
v\,(|N_h|)&=\escpr{D_{v}N, \nu_{h}}-\escpr{N,T}\,
\escpr{J(v),\nu_{h}},
\\
\label{eq:vnt}
v(\escpr{N,T})&=\escpr{D_{v}N,T}+\escpr{N,J(v)}, \\
\label{eq:dvnuh}
D_{v}\nu_h&=|N_h|^{-1}\, \big(\escpr{D_vN,Z}-\escpr{N,T}\,
\escpr{J(v),Z}\big)\,Z+\escpr{Z,v}\,T.
\end{align}
\end{lemma}

\begin{proof}
Equalities \eqref{eq:dvnh} and \eqref{eq:vmnh} are easily obtained
since $N_h=N-~\escpr{N,T}\,T$.  The proof of \eqref{eq:vnt} is
immediate.  Let us show that \eqref{eq:dvnuh} holds.  As $|\nu_h|=1$
and $\{Z_{p},(\nu_{h})_{p},T_{p}\}$ is an orthonormal basis of
$T_p\hh^1$, we get
\[
D_{v}\nu_{h}=\escpr{D_{v}\nu_{h},Z}\,Z+\escpr{D_{v}\nu_{h},T}\,T.
\]
Note that $\escpr{D_{v}\nu_{h},T}=-\escpr{\nu_{h},J(v)}=\escpr{Z,v}$ by
\eqref{eq:conmute}.  On the other hand, by using \eqref{eq:dvnh} and
the fact that $Z$ is tangent and horizontal, we deduce
\[
\escpr{D_{v}\nu_{h},Z}=|N_h|^{-1}\escpr{D_{v}N_h,Z}=
|N_h|^{-1}\big(\escpr{D_vN,Z}-\escpr{N,T}\escpr{J(v),Z}\big),
\]
and the proof follows.
\end{proof}

\begin{remark}
In $\Sg-\Sg_{0}$ we can consider the orthonormal basis $\{Z,S\}$
defined in \eqref{eq:nuh} and \eqref{eq:ese}.  By using the definition
of mean curvature in \eqref{eq:mc} we have
\[
-2H=\divv_{\Sg}\nuh=\escpr{D_{Z}\nuh,Z}+\escpr{D_{S}\nuh,S}.
\]
By \eqref{eq:dvnuh} we get $D_{Z}\nuh=T-\mnh^{-1}\escpr{B(Z),Z}Z$, and
that $D_{S}\nuh$ is proportional to $Z$.  It follows that, in
$\Sg-\Sg_{0}$
\begin{align}
\label{eq:mc2}
2H&=|N_{h}|^{-1}\,\escpr{B(Z),Z},
\\
\label{eq:dznuh}
D_{Z}\nuh&=T-(2H)\,Z,
\end{align}
where $B$ is the Riemannian shape operator of $\Sg$.  On the
other hand, the vector $D_{Z}Z$ is orthogonal to $Z$ and $T$ since
$|Z|=1$ and $\escpr{J(Z),Z}=0$.  It follows that $D_{Z}Z$ is
proportional to $\nuh$.  From \eqref{eq:dznuh} we obtain
\begin{equation}
\label{eq:dzz}
D_{Z}Z=\escpr{D_{Z}Z,\nuh}\,\nuh=2H\,\nuh.
\end{equation}
\end{remark}

The second derivative of the area for non-singular variations of a
minimal surface in $\hh^1$ has appeared in several contexts, see
\cite[Prop.~6.1]{chmy}, \cite[Sect.~3.2]{bscv}, \cite[Sect.~14]{dgn},
\cite[Proof of Thm.~3.5]{mscv} and \cite[Thm.~E]{hp2}.  In the next
theorem we compute the second derivative of the area functional for
some \emph{non-singular variations by Riemannian geodesics} of
a $C^2$ minimal surface (maybe with non-empty boundary) in $\hh^1$.

\begin{theorem}
\label{th:2ndvar}
Let $\Sg\sub\hh^1$ be a $C^2$ immersed minimal
surface with boundary $\ptl\Sg$ and singular set $\Sg_{0}$.
Consider the $C^1$ vector field $U=vN+wT$, where $N$ is
a unit normal vector to $\Sg$ and $v,\,w \in C^1_{0}(\Sg-\Sg_{0})$.
If $u=\escpr{U,N}$, then the second derivative of the area for
the variation induced by $U$ is given by
\begin{align}
\label{eq:gen2ndvar} A''(0)=&\int_{\Sg}|N_{h}|^{-1}\left\{Z(u)^2-
\big(|B(Z)+S|^2-4|N_{h}|^2\big)\,u^2\right\}d\Sg\\
\nonumber &+\int_{\Sg}\divv_{\Sg}(\xi Z)\,d\Sg
+\int_{\Sg}\divv_{\Sg}(\mu Z)\,d\Sg.
\end{align}
Here $\{Z,S\}$ is the orthonormal basis in \eqref{eq:nuh} and
\eqref{eq:ese}, $B$ is the Riemannian shape operator of $\Sg$, the
functions $\xi$ and $\mu$ are defined by
\begin{align}
\label{eq:phi}
\xi &=\escpr{N,T}\,(1-\escpr{B(Z),S})\,u^2,\\
\label{eq:mu}
\mu&=|N_h|^2\,\left(\escpr{N,T}\,(1-\escpr{B(Z),S})\,w^2
-2\escpr{B(Z),S}\,vw\right),
\end{align}
and the divergence terms are understood in the sense of
\eqref{eq:divv}.

In particular, if $\ptl\Sg$ is empty, then
\begin{equation}
\label{eq:2ndvarnosing}
A''(0)=\int_{\Sg}|N_{h}|^{-1}\left\{Z(u)^2-
\big(|B(Z)+S|^2-4|N_{h}|^2\big)\,u^2\right\}d\Sg.
\end{equation}
\end{theorem}

\begin{proof}
We will follow closely the arguments in \cite[\S
9]{simon}. Let $\var_{s}(p)=\exp_{p}(sU_{p})$, for $s$ small, be the variation induced by $U$. Then any $\Sg_{s}=\varphi_{s}(\Sg)$ is a $C^1$ immersed
oriented surface. We extend the vector $U$ along the variation by
setting $U(\varphi_{s}(p))=(d/dt)|_{t=s}\,\varphi_{t}(p)$.  Let $N$ be a continuous vector field along the
variation whose restriction to any $\Sg_{s}$ is a unit normal vector.
By using \eqref{eq:area}, the coarea formula, and that the Riemannian
area of $\Sg_{0}$ vanishes, we have
\begin{equation}
\label{eq:arform}
A(s)=A(\Sg_{s})=\int_{\Sg_{s}}\mnh\,d\Sg_{s}=
\int_{\Sg-\Sg_{0}}(\mnh\circ\var_{s}) \,|\text{Jac}\,\var_{s}|\,d\Sg,
\end{equation}
where $\text{Jac}\,\var_{s}$ is the Jacobian determinant of the
diffeomorphism $\var_{s}:\Sg\to\Sg_{s}$.

We can suppose that $\mnh(\var_s(p))>0$ whenever
$p\in\Sg-\Sg_0$ and $|s|<s_0$. Take a point $p\in\Sg-\Sg_{0}$ and consider the orthonormal basis
$\{e_1,e_{2}\}$ of $T_{p}\Sg$ given by $e_{1}=Z_{p}$ and
$e_{2}=S_{p}$.  Let $\ga$ be the Riemannian geodesic defined by
$\ga(s)=\var_{s}(p)=\exp_p(sU_p)$.  Denote by
$N(s)$ the unit normal to $\Sg_{s}$ at $\ga(s)$.  Let
$\alpha_{i}:(-\eps_{0},\eps_{0})\to\Sg-\Sg_{0}$ be a $C^1$ curve
such that $\alpha_{i}(0)=p$ and $\dot{\alpha_{i}}(0)=e_{i}$.  We
define the $C^1$ map $F_{i}:(-\eps_{0},\eps_{0})\times\rr\to\hh^1$
given by $F_{i}(\eps,s)=\var_{s}(\alpha_{i}(\eps))
=\exp_{\alpha_{i}(\eps)}(sU_{\alpha_{i}(\eps)})$.  By using
Lemma~\ref{lem:jacobi} we deduce that $E_{i}(s)=(\ptl
F_{i}/\ptl\eps)(0,s)=e_{i}(\var_{s})$ is a $C^\infty$ Jacobi
vector field along $\ga$ with $[\dot{\ga},E_{i}]=0$ and
$E_{i}(0)=e_{i}$. Therefore, we have the following identities
along $\ga$
\begin{align}
\label{eq:jacobi} D_{U}D_{U}E_{i}&=-R(U,E_{i})U,
\\
\label{eq:bracket} D_{U}E_{i}&=D_{E_{i}}U.
\end{align}
On the other hand, it is clear that $\{E_{1}(s),E_{2}(s)\} $
provide a basis of the tangent space to $\Sg_{s}$ at $\ga(s)$.  In
particular $|\text{Jac}\,\var_{s}|=
(|E_{1}|^2\,|E_{2}|^2-\escpr{E_{1},E_{2}}^2)^{1/2}(s)$, and so
$|\text{Jac}\,\var_{s}|$ is $C^\infty$ along $\ga$.  Moreover, we
have $N(s)=\pm|E_{1}\times E_{2}|^{-1}\,(E_{1}\times E_{2})(s)$,
which is $C^\infty$ on $\ga$.  Here $\times$ is the cross product
in $(\hh^1,g)$.  We conclude that $|N_{h}|(s)$ is $C^\infty$ along
$\ga$ as well. Thus we can apply the classical result of
differentiation under the integral sign to deduce, from
\eqref{eq:arform}, that
\begin{equation}
\label{eq:2nd1}
A''(0)=\int_{\Sg-\Sg_{0}}\left\{\mnh''(0)
+2\,\mnh'(0)\,|\text{Jac}\,\var_{s}|'(0)
+\mnh\,|\text{Jac}\,\var_{s}|''(0)\right\}d\Sg,
\end{equation}
where we have used that $\var_{0}(p)=p$ for any $p\in\Sg$, and so
$|\text{Jac}\,\var_{0}|=1$.

Now we compute the different terms in \eqref{eq:2nd1}.  The
calculus of $|\text{Jac}\,\var_{s}|'(0)$ and
$|\text{Jac}\,\var_{s}|''(0)$ is found in \cite[\S 9]{simon} for
$C^2$ variations of a $C^1$ surface in Euclidean space.  The
arguments can be generalized to any Riemannian manifold for a
$C^1$ variation obtained when we leave from a $C^2$ surface by
geodesics.  As $U=vN+wT$ on $\Sg$, we deduce, by using
$\divv_\Sg T=0$ and the third equality in \eqref{eq:relations}, that
\begin{equation}
\label{eq:jac1}
|\text{Jac}\,\var_{s}|'(0)=\divv_{\Sg}U=(-2H_{R})\,v-|N_h|\,S(w)
=-\escpr{B(S),S}\,v-|N_h|\,S(w).
\end{equation}
To get the second equality we have taken into account
\eqref{eq:mc2} to obtain
\[
2H_{R}=-\divv_{\Sg}N=\escpr{B(Z),Z}+\escpr{B(S),S}=
\escpr{B(S),S}.
\]
On the other hand, it is known that
\begin{align*}
|\text{Jac}\,\var_{s}|''(0)=&(\divv_{\Sg}U)^2
+\sum_{i=1}^2|(D_{e_{i}}U)^\bot|^2
\\
\nonumber
& -\sum_{i=1}^2\escpr{R(U,e_{i})U,e_{i}}
-\sum_{i,j=1}^2\escpr{D_{e_{i}}U,e_{j}}\,\escpr{D_{e_{j}}U,e_{i}}.
\end{align*}
Hence from \eqref{eq:jac1}, equality
\begin{equation}
\label{eq:deu}
D_eU=e(v)\,N-v\,B(e)+e(w)\,T+w\,J(e),
\end{equation}
and equation \eqref{eq:curvature}, we get
\begin{align}
\label{eq:jac2}
|\text{Jac}\,\var_{s}|''(0)=&|\nabla_{\Sg}v|^2
+\escpr{N,T}^2\,|\nabla_{\Sg}w|^2
-2|N_h|\,Z(v)\,w-2|N_h|\,\escpr{B(Z),S}\,Z(w)\,v
\\
\nonumber
&+2\escpr{N,T}\,Z(v)\,Z(w)+2\escpr{N,T}\,S(v)\,S(w)
\\
\nonumber
&-(\ric(N,N)+|B|^2-\escpr{B(S),S}^2)\,v^2-4\escpr{N,T}\,vw,
\end{align}
where $\nabla_{\Sg}$ is the gradient relative
to $\Sg$, $\ric$ is the Ricci tensor in $(\hh^1,g)$, and
$|B|^2$ is the squared norm of the Riemannian shape operator of
$\Sg$.

Let us compute $\mnh'(0)$ and $\mnh''(0)$.  From \eqref{eq:vmnh}
and \eqref{eq:conmute} it follows that
\[
\mnh'(s)=U(\mnh)=\escpr{D_{U}N,\nuh}
-\escpr{N,T}\,\escpr{J(U),\nuh}=\escpr{D_{U}N,\nuh}
+\escpr{N,T}\,\escpr{U,Z}.
\]
Note that $U=uN-(\mnh w)S$. Then $\escpr{U,Z}=0$ and
$D_UN=-\nabla_\Sg u+(\mnh w)B(S)$ on $\Sg-\Sg_0$. By the second
equality in \eqref{eq:relations} we obtain
\begin{equation}
\label{eq:dmnh}
|N_{h}|'(0)=\escpr{D_{U}N,\nuh}=-\escpr{N,T}\,S(u)+|N_h|\,\escpr{N,T}\,\escpr{B(S),S}\,w.
\end{equation}
We also deduce the following
\begin{align}
\label{eq:d2mnh1}
|N_{h}|''(0)=&\escpr{D_{U}D_{U}N,\nuh}+\escpr{D_{U}N,D_{U}\nuh}
\\
\nonumber
&+U(\escpr{N,T})\,\escpr{U,Z}+\escpr{N,T}\,U(\escpr{U,Z})
\\
\nonumber =&\escpr{D_{U}D_{U}N,\nuh}+\escpr{D_{U}N,D_{U}\nuh}
+\escpr{N,T}\,\escpr{U,D_{U}Z},
\end{align}
since $\escpr{U,Z}=0$ and $D_{U}U=0$ on $\Sg-\Sg_{0}$.  We can
compute $D_{U}\nuh$ from \eqref{eq:dvnuh}.  By using that
$D_{U}N=-\nabla_{\Sg}u+(|N_h|w)B(S)$ and $J(U)=(\mnh v)Z$ on
$\Sg-\Sg_{0}$, we get
\[
D_{U}\nuh=-|N_{h}|^{-1}\,(Z(u)-|N_h|\,\escpr{B(Z),S}\,w+|N_{h}|\,\escpr{N,T}\,v)\,Z,
\]
and so
\begin{align}
\label{eq:d2mnh12}
\escpr{D_{U}N,D_{U}\nuh}=&\mnh^{-1}\,Z(u)^2+\escpr{N,T}\,Z(u)\,v
-2\escpr{B(Z),S}\,Z(u)\,w\\
\nonumber
&-|N_h|\,\escpr{N,T}\,\escpr{B(Z),S}\,vw+|N_h|\,\escpr{B(Z),S}^2\,w^2.
\end{align}
Now we compute $D_{U}Z$.  The coordinates of this vector with
respect to the orthonormal basis $\{Z,\nuh,T\}$ are given by
\[
\escpr{D_{U}Z,Z}=0, \quad
\escpr{D_{U}Z,\nuh}=|N_{h}|^{-1}\,Z(u)-\escpr{B(Z),S}\,w+
\escpr{N,T}\,v,\quad\escpr{D_{U}Z,T}=-|N_{h}|\,v.
\]
The previous equalities and the fact that $U=vN+wT$ on
$\Sg-\Sg_{0}$ imply that
\begin{equation}
\label{d2mnh13}
\escpr{U,D_{U}Z}=Z(u)\,v-\mnh\,(1+\escpr{B(Z),S})\,vw.
\end{equation}
It remains to compute $D_{U}D_{U}N$.  Note that
$\{E_{1},E_{2},N\}$ provides an orthonormal basis of $T_{p}\hh^1$.
As a consequence
\[
D_{U}D_{U}N=\sum_{i=1}^2\escpr{D_{U}D_{U}N,E_{i}}\,E_{i}
+\escpr{D_{U}D_{U}N,N}\,N.
\]
As $\escpr{N,E_{i}}=0$ along $\ga$  we get
\begin{align*}
\escpr{D_{U}D_{U}N,E_{i}}&=-2\escpr{D_{U}N,D_{U}E_{i}}
-\escpr{N,D_{U}D_{U}E_{i}}
\\
\nonumber
&=-2\escpr{D_{U}N,D_{E_{i}}U}+\escpr{N,R(U,E_{i})U}.
\end{align*}
The second equality follows from \eqref{eq:bracket} and
\eqref{eq:jacobi}. Now recall that $e_{1}=Z_{p}$ and
$e_{2}=S_{p}$. It follows that
\begin{align}
\label{eq:d2mnh11}
\escpr{D_{U}D_{U}N,\nuh}=&-2
\escpr{N,T}\,\escpr{D_{U}N,D_{S}U} +
\escpr{N,T}\,\escpr{N,R(U,S)U}
\\
\nonumber
&+|N_h|\,\escpr{D_{U}D_{U}N,N}.
\end{align}
By taking into account \eqref{eq:deu} and that
$D_UN=-\nabla_{\Sg}u+(|N_h|\,w)B(S)$, we obtain
\begin{align}
\label{eq:dundsu}
\escpr{D_{U}N,D_{S}U}=&
-|N_h|^2\,\escpr{B(S),S}\,S(w)\,w-\escpr{N,T}\,Z(u)\,w
\\
\nonumber
&+|N_h|\,S(u)\,S(w)+B(S)(u)\,v
\\
\nonumber
&+|N_h|\,\escpr{N,T}\,\escpr{B(Z),S}\,w^2-|N_h|\,|B(S)|^2\,vw.
\end{align}
On the other hand, we use \eqref{eq:curvature} so that, after a
straightforward computation, we conclude
\begin{equation}
\label{eq:rusun}
\escpr{R(U,S)U,N}=|N_h|\,(v+\escpr{N,T}\,w)\,w.
\end{equation}
Moreover, since $|N|^2=1$ on $\Sg-\Sg_{0}$ we have
\begin{equation}
\label{eq:d2n12}
\escpr{D_{U}D_{U}N,N}=-|D_{U}N|^2=
-|\nabla_{\Sg}u|^2+2|N_h|\,B(S)(u)\,w-|N_h|^2\,|B(S)|^2\,w^2.
\end{equation}

By substituting \eqref{eq:dundsu}, \eqref{eq:rusun} and \eqref{eq:d2n12}
into \eqref{eq:d2mnh11} we get $\escpr{D_UD_UN,\nuh}$.
From \eqref{eq:d2mnh11}, \eqref{eq:d2mnh12} and \eqref{d2mnh13},
after simplifying, equality \eqref{eq:d2mnh1} becomes
\begin{align}
\label{eq:d2mnh2}
\mnh''(0)=&\mnh^{-1}\,Z(u)^2-\mnh\,|\nabla_{\Sg}u|^2
+2\escpr{N,T}\,(Z(u)\,v-B(S)(u)\,v)
\\
\nonumber
&+2\big(\escpr{N,T}^2+\mnh^2\,\escpr{B(Z),S}
-\escpr{B(Z),S}\big)\,Z(u)\,w
\\
\nonumber
&+2\mnh^2\,\escpr{B(S),S}\,\big(S(u)\,w+\escpr{N,T}\,S(w)\,w\big)
-2\mnh\,\escpr{N,T}\,S(u)\,S(w)
\\
\nonumber
&+2|N_h|\,\escpr{N,T}\,\big(|B(S)|^2-\escpr{B(Z),S}\big)\,vw
\\
\nonumber
&+\big(|N_h|\,\escpr{N,T}^2\,(1-\escpr{B(Z),S})^2
-|N_h|^3\,\escpr{B(S),S}^2\big)\,w^2.
\end{align}

Now, since $u=v+\escpr{N,T}w$, we have $\nabla_\Sg u=\nabla_\Sg v
+w\,\nabla_\Sg(\escpr{N,T})+\escpr{N,T}\,\nabla_\Sg w$. By \eqref{eq:vnt}
and \eqref{eq:relations} it is easy to see that
\begin{align}
\label{eq:znt}
Z(\escpr{N,T})&=\mnh\,(\escpr{B(Z),S}-1),
\\
\nonumber
S(\escpr{N,T})&=\mnh\,\escpr{B(S),S}.
\end{align}
This allows us to compute the term $\mnh\,|\nabla_\Sg u|^2$ in
\eqref{eq:d2mnh2}.
At this moment, we use \eqref{eq:d2mnh2}, \eqref{eq:dmnh},
\eqref{eq:jac1} and \eqref{eq:jac2} so that, after simplifying, we get that
\begin{equation*}
\mnh''(0)
+2\,\mnh'(0)\,|\text{Jac}\,\var_{s}|'(0)
+\mnh\,|\text{Jac}\,\var_{s}|''(0)
\end{equation*}
is equal to
\begin{align}
\label{eq:2nd2}
\mnh^{-1}&\,Z(u)^2
+2\escpr{N,T}\, (1-\escpr{B(Z),S}) \,(Z(v)\,v+Z(w)\,w)\\
\nonumber &+2(\escpr{N,T}^2-\escpr{B(Z),S})\,(Z(w)\,v+Z(v)\,w)
+q_{1}\,v^2
\\
\nonumber
& +2 |N_h|\,\escpr{N,T}\,(\escpr{B(Z),S}-3)\,vw
-|N_h|(1-\escpr{B(Z),S})^2\,w^2,
\end{align}
where $q_{1}$ is the function given by
\[
q_{1}=|N_{h}|\,\big(\escpr{B(S),S}^2-\ric(N,N)-|B|^2\big).
\]

In order to obtain \eqref{eq:gen2ndvar} from \eqref{eq:2nd1} and
\eqref{eq:2nd2}, we apply Lemma~\ref{lem:aux4} below.  We
deduce the following
\begin{align*}
\mnh''(0)
&+2\,\mnh'(0)\,|\text{Jac}\,\var_{s}|'(0)
+\mnh\,|\text{Jac}\,\var_{s}|''(0)
\\
=&\mnh^{-1}\,Z(u)^2+\divv_{\Sg}(\rho Z)
\\
\nonumber
&+\left\{q_{1}+\big(\escpr{B(Z),S}-1\big)
\big(\escpr{N,T}\,q_2+Z(\escpr{N,T})\big)
+\escpr{N,T}\,Z(\escpr{B(Z),S})
\right\}v^2
\\
\nonumber
&+\left\{\big(\escpr{B(Z),S}-1\big)
\big(\escpr{N,T}\,q_2+Z(\escpr{N,T})\big)
+\escpr{N,T}\,Z(\escpr{B(Z),S})
\right.
\\
&\left.
-\mnh\,(1-\escpr{B(Z),S})^2
\right\}w^2+2\left\{\mnh\,\escpr{N,T}\,(\escpr{B(Z),S}-3)
-Z(\escpr{N,T}^2)\right.
\\
\nonumber
&\left.-\escpr{N,T}^2\,q_2
+\escpr{B(Z),S}\,q_2+Z(\escpr{B(Z),S})
\right\} vw,
\end{align*}
where $\rho$ is the function
\[
\escpr{N,T}\,(1-\escpr{B(Z),S})\,(v^2+w^2)
+2\,(\escpr{N,T}^2-\escpr{B(Z),S})\,vw.
\]
A straightforward computation using \eqref{eq:znt}, \eqref{eq:zbzs},
the identities
\begin{align}
\nonumber
\ric(N,N)&=2-4\mnh^2\quad\text{(it follows from \eqref{eq:ricci})},
\\
\nonumber
|B|^2&=\escpr{B(Z),Z}^2+\escpr{B(S),S}^2+2\escpr{B(Z),S}^2=
\escpr{B(S),S}^2+2\escpr{B(Z),S}^2,
\\
\nonumber
B(Z)&=\escpr{B(Z),Z}\,Z+\escpr{B(Z),S}\,S=\escpr{B(Z),S}\,S,
\end{align}
and that $u=v+\escpr{N,T}w$, gives us
\begin{align*}
\mnh''(0) +2\,\mnh'(0)\,&|\text{Jac}\,\var_{s}|'(0)
+\mnh\,|\text{Jac}\,\var_{s}|''(0)
\\
 =&\mnh^{-1}\,Z(u)^2-|N_{h}|^{-1}\,(|B(Z)+S|^2-4|N_{h}|^2)\,u^2\\
 &+\divv_{\Sg}(\xi Z)+\divv_{\Sg}(\mu Z),
\end{align*}
where $\xi$ and $\mu$ are the functions given in \eqref{eq:phi}
and \eqref{eq:mu}.

Finally, suppose that $\ptl\Sg$ is empty. Then
$\xi$ and $\mu$ are continuous functions with compact support in
$\Sg-\Sg_{0}$ and they are also $C^1$ in the $Z$-direction by
Lemma~\ref{lem:zbzs}. Hence the integrals of $\divv_{\Sg}(\xi Z)$
and  $\divv_{\Sg}(\mu Z)$ vanish by virtue of the divergence
theorem in Lemma~\ref{lem:div}.  This proves
\eqref{eq:2ndvarnosing}.
\end{proof}

\begin{remark}
The divergence terms in \eqref{eq:gen2ndvar} need not vanish
if $\ptl\Sg$ is nonempty. In the proof of
Proposition~\ref{prop:2ndvarsing} we will show that these terms
play an important role.
\end{remark}

\begin{lemma}
\label{lem:aux4}
Let $\Sg$ be a $C^2$ immersed oriented surface in
$\hh^1$ and $\phi\in C^1(\Sg)$.  Then, in the regular set
$\Sg-\Sg_{0}$, we have
\[
\divv_{\Sg}(\phi\,Z)=Z(\phi)+q_{2}\,\phi,
\]
where $q_{2}$ is the function given by
\[
q_{2}=\mnh^{-1}\escpr{N,T}\,(1+\escpr{B(Z),S}).
\]
\end{lemma}

\begin{proof}
Clearly we have
\begin{equation}
\label{eq:firstoso}
\divv_{\Sg}(\phi\,Z)=(\divv_{\Sg}Z)\,\phi+Z(\phi).
\end{equation}
Note that
\begin{equation*}
\label{eq:toto21}
\divv_{\Sg}Z=\escpr{D_{Z}Z,Z}+\escpr{D_{S}Z,S}=
\escpr{D_{S}Z,S},
\end{equation*}
since $|Z|^2=1$.  We compute the components of $D_{S}Z$ in the
orthonormal basis $\{Z,\nuh,T\}$.  Observe that $D_{S}Z$ is orthogonal to
$Z$.  By using \eqref{eq:dvnuh} and that $J(S)=\escpr{N,T}Z$, we get
\begin{align*}
\escpr{D_{S}Z,\nuh}&=-\escpr{Z,D_{S}\nuh}
=|N_{h}|^{-1}\,(\escpr{B(Z),S}+\escpr{N,T}^2),
\\
\escpr{D_{S}Z,T}&=-\escpr{Z,J(S)}=-\escpr{N,T}.
\end{align*}
From here we deduce
\begin{equation}
\label{eq:dsz}
D_{S}Z=\mnh^{-1}\,(\escpr{B(Z),S}+1-|N_{h}|^2)\,\nuh-\escpr{N,T}\,T.
\end{equation}
As a consequence, we obtain
\begin{equation}
\label{eq:toto12}
\divv_{\Sg}Z=\escpr{D_{S}Z,S}=
\mnh^{-1}\escpr{N,T}\,(1+\escpr{B(Z),S}).
\end{equation}
The proof finishes by substituting \eqref{eq:toto12}
into \eqref{eq:firstoso}.
\end{proof}

\begin{lemma}
\label{lem:zbzs}
Let $\Sg$ be a $C^2$ immersed oriented minimal surface in $\hh^1$.
Then, in the regular set $\Sg-\Sg_{0}$, we have
\begin{itemize}
\item[(i)] The functions $\escpr{N,T}$ and $\mnh$ are
$C^\infty$ in the $Z$-direction.
\item[(ii)] The vector fields $\nuh$ and $S$ are $C^\infty$
in the $Z$-direction.
\item[(iii)] The function $\escpr{B(Z),S}$ is $C^\infty$
in the $Z$-direction, and
\begin{equation}
\label{eq:zbzs}
Z(\escpr{B(Z),S})=4|N_{h}|\,\escpr{N,T}-
2|N_{h}|^{-1}\escpr{N,T}\,\escpr{B(Z),S}\,(1+\escpr{B(Z),S}).
\end{equation}
\end{itemize}
\end{lemma}

\begin{proof}
Recall that $N$ is $C^\infty$ in the $Z$-direction by
Lemma~\ref{lem:reg}.  This implies (i).  Assertions (ii) and
(iii) follow from (i) by the definition of $\nuh$ and $S$ in
\eqref{eq:nuh} and \eqref{eq:ese}.  To compute $Z(\escpr{B(Z),S})$
note that
\[
Z(\escpr{B(Z),S})=Z(-\escpr{D_{Z}N,S})=
-\escpr{D_{Z}D_{Z}N,S}-\escpr{D_{Z}N,D_{Z}S}.
\]
It is clear that $D_{Z}N$ is tangent to $\Sg$.  On the other hand,
$D_{Z}S$ is proportional to $N$.  This comes from the fact that
$\escpr{D_{Z}S,Z}=-\escpr{S,D_{Z}Z}=0$ by \eqref{eq:dzz}, whereas
$\escpr{D_{Z}S,S}=0$.  Therefore we have
\begin{align}
\label{eq:dzndzs}
\escpr{D_{Z}N,D_{Z}S}&=0,
\\
\label{eq:zbzs1}
Z(\escpr{B(Z),S})&=-\escpr{D_{Z}D_{Z}N,S}=\escpr{N,D_{Z}D_{Z}S}.
\end{align}

It remains to compute $D_{Z}D_{Z}S$.  From \eqref{eq:dsz} we see that
$D_{S}Z$ is $C^\infty$ in the $Z$-direction.  As a consequence
$[Z,S]=D_{Z}S-D_{S}Z$ is also $C^\infty$ in the $Z$ direction, and
$D_{Z}[Z,S]=D_{Z}D_{Z}S-D_{Z}D_{S}Z$.  Thus equation \eqref{eq:zbzs1}
becomes
\begin{align}
\label{eq:zbzs2}
Z(\escpr{B(Z),S})&=\escpr{N,D_{Z}[Z,S]}+\escpr{N,D_{Z}D_{S}Z}
\\
\nonumber
&=\escpr{N,D_{Z}[Z,S]}+\escpr{N,D_{S}D_{Z}Z}-\escpr{N,R(Z,S)Z}
+\escpr{N,D_{[Z,S]}Z}
\\
\nonumber &=\escpr{N,D_{Z}[Z,S]}-\escpr{N,R(Z,S)Z}
+\escpr{N,D_{[Z,S]}Z},
\end{align}
where $R$ is the Riemannian curvature tensor and we have used
\eqref{eq:dzz} to get $D_{S}D_{Z}Z=0$.  Now, observe that
\[
\escpr{[Z,S],N}=\escpr{D_{Z}S,N}-\escpr{D_{S}Z,N}=-\escpr{S,D_{Z}N}+
\escpr{Z,D_{S}N}=0,
\]
which implies that $[Z,S]$ is tangent to $\Sg$.  Therefore, we deduce
\begin{align*}
\escpr{N,D_{Z}[Z,S]}&=\escpr{B(Z),[Z,S]}=\escpr{B(Z),D_{Z}S}
-\escpr{B(Z),D_{S}Z}=-\escpr{B(Z),D_{S}Z},
\\
\escpr{N,D_{[Z,S]}Z}&=-\escpr{D_{[Z,S]}N,Z}=\escpr{B(Z),[Z,S]}=
-\escpr{B(Z),D_{S}Z},
\end{align*}
where we have used \eqref{eq:dzndzs}.  If we put this information into
\eqref{eq:zbzs2}, we obtain
\begin{equation}
\label{eq:zbzs3}
Z(\escpr{B(Z),S})=-2\escpr{B(Z),D_{S}Z}-\escpr{N,R(Z,S)Z}.
\end{equation}
To compute the first term above we take into account \eqref{eq:dsz}.
After simplifying, we get
\begin{equation}
\label{eq:bzdsz}
\escpr{B(Z),D_{S}Z}=\mnh^{-1}\,\escpr{N,T}\,\escpr{B(Z),S}\,
(1+\escpr{B(Z),S}).
\end{equation}
For the second term, we apply \eqref{eq:curvature} so that, after a
straightforward calculus, we conclude
\begin{equation}
\label{eq:rzszn}
\escpr{N,R(Z,S)Z}=-4\mnh\,\escpr{N,T}.
\end{equation}
The proof finishes by substituting \eqref{eq:bzdsz} and
\eqref{eq:rzszn} into \eqref{eq:zbzs3}.
\end{proof}

In the next result we compute the second derivative of the area for some
\emph{vertical variations} of an area-stationary surface $\Sg$ whose singular
cuves $(\Sg_0)_c$ are $C^3$ (in \cite[Prop.~4.20]{rr2} we proved that they are always $C^2$). We suppose that the variation is constant along the
characteristic curves of a tubular neighborhood around $(\Sg_0)_c$. By a tubular neighborhood
of radius $\eps>0$ we mean the union of all the characteristic segments of
length $2\eps$ centered at $(\Sg_0)_c$.

\begin{proposition}
\label{prop:2ndvert}
Let $\Sg$ be a $C^2$ immersed oriented area-stationary surface
in $\hh^1$ such that the singular curves $(\Sg_0)_c$ of $\Sg$ are
of class $C^3$. Let $\var_r(p):=\exp_p(rw(p)T_p)$, for $r$ small, be
the vertical variation of $\Sg$ induced by a function $w\in C^2_0(\Sg)$.
Suppose that there is a tubular neighborhood $E_0$ of
$\supp(w)\cap (\Sg_0)_c$ where $Z(w)=0$. Then, there is
a tubular neighborhood $E$ of $\supp(w)\cap (\Sg_0)_c$ such that
\[
\frac{d^2}{dr^2}\bigg|_{r=0} A(\varphi_r(E))=\int_{(\Sg_0)_c}
S(w)^2\,dl,
\]
where $S$ is any continuous extension of the vector field $S$
defined in \eqref{eq:ese} to $(\Sg_0)_c$ and $dl$ denotes the Riemannian
length element.
\end{proposition}

\begin{proof}
We can restrict ourselves to a neighborhood of a single singular
curve $\Ga$. We consider a parameterization $\Ga(\eps)=(x(\eps),y(\eps),t(\eps))$ by arc-length. By Proposition~\ref{prop:varprop} the area-stationary surface $\Sg$
can be parameterized in a neighborhood of $\supp(w)\cap\Ga$ by
\[
(\eps,s)\mapsto \Ga(\eps)+sJ(\dot{\Ga}(\eps)),
\]
so that the curves with $\eps$ constant are the characteristic
curves of $\Sg$. In Euclidean coordinates we have
\begin{align*}
x(\eps,s)&=x(\eps)-s\dot{y}(\eps), \\
y(\eps,s)&=y(\eps)+s\dot{x}(\eps), \\
t(\eps,s)&=t(\eps)-s\,(x\dot{x}+y\dot{y})(\eps).
\end{align*}
As $Z(w)=0$ we get that $w$ is a function of $\eps$ alone. The deformation
$\var_r(p)=\exp_p(rw(p)T_p)$ consists on changing the
$t$-coordinate of the above parameterization by
\[
t(\eps,s)+rw(\eps).
\]
A simple computation shows that the tangent space to the surface
$\Sg_r:=\var_r(\Sg)$ is generated by the vectors
\begin{equation}
\label{eq:spansgr} -\dot{y}X+\dot{x}Y, \qquad
(\dot{x}-s\ddot{y})\,X+(\dot{y}+s\ddot{x})\,Y+\big(s(-2+sh)+r\dot{w}\big)\,T,
\end{equation}
where $x$, $y$ and $t$ are the coordinates of $\Ga$, dots
represent derivatives with respect to $\eps$, and
$h=h(\eps)=(\dot{x}\ddot{y}-\dot{y}\ddot{x})(\eps)$ is the
Euclidean geodesic curvature of the $xy$-projection of $\Ga$.

Hence the singular points of $\Sg_r$ corresponds to the zero set
of $F(\eps,s,r):=s(-2+sh(\eps))+r\dot{w}(\eps)$. Observe that $F$
is a $C^1$ function since the singular curves are assumed to be of
class $C^3$ and $w\in C^2$. As $(\ptl F/\ptl s)(\eps,0,0)=-2$, we can apply the
Implicit Function Theorem and a compactness argument to show that
there are positive values $\eps_0$, $s_0$, $r_0$, and a $C^1$ function
$s:(-\eps_0,\eps_0)\times (-r_0,r_0)\to (-s_0,s_0)$ with $s(\eps,0)=0$
satisfying $F(\eps,s(\eps,r),r)=0$. Here $\eps_0>0$ is taken so
that $\supp(w)\cap\Ga\subset [-\eps_0,\eps_0]$. We define
$E:=F((-\eps_0,\eps_0)\times (-s_0,s_0))$.

On the other hand, a computation using \eqref{eq:spansgr} shows
that
\[
|(N_h)_r|\,d\Sg_r=|s(-2+sh(\eps))+r\dot{w}(\eps)|\,d\eps\,ds.
\]
Hence we have
\[
A(\varphi_r(E))=\int_{-\eps_0}^{\eps_0}\bigg\{\int_{-s_0}^{s_0}
|s(-2+sh(\eps))+r\dot{w}(\eps)|\,ds\bigg\}\,d\eps.
\]
Denote by $f_\eps(r)$ the integral between brackets.
As $(\ptl F/\ptl s)(\eps,0,0)<0$ we deduce
\[
f_\eps(r)=\int_{-s_0}^{s(\eps,r)}
\big(s(-2+sh(\eps))+r\dot{w}(\eps)\big)\,ds +
\int_{s(\eps,r)}^{s_0} \big(s(2-sh(\eps))-r\dot{w}(\eps)\big)\,ds.
\]
Taking derivatives with respect to $r$ we obtain
\[
f_\eps'(r)=\int_{-s_0}^{s(\eps,r)}\dot{w}(\eps)\,ds
-\int_{s(\eps,r)}^{s_0}\dot{w}(\eps)\,ds=
2\dot{w}(\eps)s(\eps,r).
\]
Taking derivatives again we have
\[
f_\eps''(r)=2\dot{w}(\eps)\,\frac{\ptl s}{\ptl r}(\eps,r).
\]
Since $(\ptl s/\ptl r)(\eps,0)=\dot{w}(\eps)/2$ we conclude
\[
f_\eps''(0)=\dot{w}(\eps)^2,
\]
and so
\[
\frac{d^2}{dr^2}\bigg|_{r=0}
A(\varphi_r(E))=\int_{-\eps_0}^{\eps_0} \dot{w}(\eps)^2\,d\eps.
\]
By Proposition~\ref{prop:varprop} we know that the vector field
$S$ defined in \eqref{eq:ese} extends continuously to $\Ga$ as a unit tangent vector to
$\Ga$. Then $\dot{w}(\eps)^2=S(w)^2$ and the claim follows.
\end{proof}

\subsection{A stability criterion for stable surfaces in  $\hh^1$}
\label{subsec:sc}
Here we obtain a useful criterion to check if a
given area-stationary surface is unstable.  First we need a
definition.  Let $\Sg$ be a $C^2$ oriented minimal surface immersed in
$\hh^1$.  For two functions $u,v\in C_{0}(\Sg-\Sg_{0})$ which are also
$C^1$ in the $Z$-direction, we denote
\begin{equation}
\label{eq:indexform}
\mathcal{I}(u,v):=\int_{\Sg}\mnh^{-1}\left\{Z(u)\,Z(v)-
\big(|B(Z)+S|^2-4\mnh^2\big)\,uv\right\}\,d\Sg,
\end{equation}
where $\{Z,S\}$ is the orthonormal basis in \eqref{eq:nuh} and
\eqref{eq:ese}, and $B$ is the Riemannian shape operator of
$\Sg$.  The expression \eqref{eq:indexform} defines a symmetric
bilinear form, which we call the \emph{index form} associated to $\Sg$
by analogy with the Riemannian situation, see \cite{bdce}.

\begin{proposition}
\label{prop:stcond1}
Let $\Sg$ be a $C^2$ immersed oriented area-stationary surface in
$\hh^1$ with singular set $\Sg_{0}$.  If $\Sg$ is stable under
non-singular variations then the index form defined in
\eqref{eq:indexform} satisfies $\mathcal{I}(u,u)\geq 0$ for any
function $u\in C_{0}(\Sg-\Sg_{0})$ which is also $C^1$ in the
direction of the characteristic field $Z$.
\end{proposition}

\begin{proof}
Let $N$ be the unit normal vector to $\Sg$.  Take $u\in
C^1_{0}(\Sg-\Sg_{0})$ and consider the vector field $U=uN$.  Note that
$\Sg$ is a minimal surface since it is area-stationary.  Hence
Theorem~\ref{th:2ndvar} implies that the second derivative of the
area for the variation induced by $U$ is $A''(0)=\mathcal{I}(u,u)$.  As
$\Sg$ is stable under non-singular variations we deduce that
\begin{equation}
\label{eq:c1}
\mathcal{I}(u,u)\geq 0,\quad\text{ for any } u\in C^1_{0}(\Sg-\Sg_{0}).
\end{equation}

Now fix a function $u\in C_{0}(\Sg-\Sg_{0})$ which is also $C^1$ in
the $Z$-direction.  By using Lemma~\ref{lem:approx} and that $\Sg_{0}$
has vanishing Riemannian area, we can find a compact set
$K\subseteq\Sg-\Sg_{0}$ and a sequence of functions
$\{u_{\eps}\}_{\eps>0}$ in $C^1_{0}(\Sg-\Sg_{0})$ such that
$\{u_{\eps}\}\to u$ in $L^2(\Sg)$, $\{Z(u_{\eps})\}\to Z(u)$ in
$L^2(\Sg)$, while the supports of $u_{\eps}$ and $u$ are contained in
$K$ for any $\eps>0$.  From here it is not difficult to check that
$\{\mnh^{-1/2}Z(u_{\eps})\}\to\mnh^{-1/2}Z(u)$, $\{(\mnh^{-1} f_{1})
^{1/2}u_{\eps}\}\to (\mnh^{-1}f_{1})^{1/2}u$ and $\{(\mnh^{-1}
f_{2})^{1/2}u_{\eps}\}\to (\mnh^{-1}f_{2})^{1/2}u$ in $L^2(\Sg)$,
where $f_{1}=|B(Z)+S|^2$ and $f_{2}=4\mnh^2$.  It follows that
$\lim_{\eps\to 0}\mathcal{I}(u_{\eps},u_{\eps})=\mathcal{I}(u,u)$, so
that inequality \eqref{eq:c1} proves the claim.
\end{proof}

\begin{remark}
As in \cite[Thm.~15.2]{dgn} and \cite[Thm.~3.5, Cor.~3.7]{mscv} the
previous result can be seen as a Poincar\'e type inequality for stable
surfaces in $\hh^1$.
\end{remark}

\subsection{Integration by parts.  The stability operator in $\hh^1$}
\label{subsec:ibp}
In Riemannian geometry the index form of a minimal surface can be
expressed in terms of a second order elliptic operator defined on the
surface, see \cite{bdce}.  In this part of the section we prove a
similar property for the index form \eqref{eq:indexform} of a minimal
surface in $\hh^1$ which involves a hypoelliptic second order
differential operator on the surface.

\begin{proposition}[Integration by parts I]
\label{prop:ibp}
Let $\Sg\sub\hh^1$ be a $C^2$ immersed surface with unit normal
vector $N$ and singular set $\Sg_{0}$.  Consider two functions $u\in
C_{0}(\Sg-\Sg_{0})$ and $v\in C(\Sg-\Sg_{0})$ which are $C^1$ and
$C^2$ in the $Z$-direction, respectively.  Then we have
\[
\mathcal{I}(u,v)=-\int_{\Sg}u\,\mathcal{L}(v)\,d\Sg,
\]
where $\mathcal{I}$ is the index form defined in \eqref{eq:indexform},
and $\mathcal{L}$ is the second order differential operator
\begin{align}
\label{eq:lu}
\mathcal{L}(v):=|N_{h}|^{-1}\big\{Z(Z(v))
&+2\,|N_{h}|^{-1}\,\escpr{N,T}\,\escpr{B(Z),S}\,Z(v)
\\
\nonumber
&+(|B(Z)+S|^2-4|N_{h}|^2)\,v\big\}.
\end{align}
\end{proposition}

\begin{proof}
Along this proof we shall denote $q=|B(Z)+S|^2-4\mnh^2$.  First note
that in $\Sg-\Sg_{0}$ the hypotheses about $u$ and $v$ ensure that
$\mnh^{-1}Z(v)$ and $\mnh^{-1}Z(v)u$ are $C^1$ in the $Z$-direction.
Suppose proved that
\begin{equation}
\label{eq:divform}
\mathcal{L}(v)=\divv_{\Sg}(\mnh^{-1}Z(v)\,Z)+\mnh^{-1}qv.
\end{equation}
In such a case, we would apply the divergence theorem in
Lemma~\ref{lem:div} in order to get
\begin{align*}
0=\int_{\Sg}\divv_{\Sg}(|N_{h}|^{-1}Z(v)\,u\,Z)\,d\Sg&=
\int_{\Sg}u\,\divv_{\Sg}(|N_{h}|^{-1}Z(v)\,Z)\,d\Sg
+\int_{\Sg}|N_{h}|^{-1}Z(u)\,Z(v)\,d\Sg.
\\
&=\int_{\Sg}u\,\mathcal{L}(v)\,d\Sg+\mathcal{I}(u,v),
\end{align*}
and this would finish the proof.

To obtain \eqref{eq:divform} observe that
\begin{equation}
\label{eq:toto2}
\divv_{\Sg}(\mnh^{-1}Z(v)\,Z)=\mnh^{-1}Z(v)\,
\divv_{\Sg}Z+Z(\mnh^{-1}Z(v)).
\end{equation}
The computation of $\divv_{\Sg}Z$ is given in \eqref{eq:toto12}.  On
the other hand, we have
\begin{align}
\label{eq:toto22}
Z(\mnh^{-1}Z(v))&=\mnh^{-1}Z(Z(v))+Z(\mnh^{-1})\,Z(v)
\\
\nonumber
&=\mnh^{-1}Z(Z(v))-\mnh^{-2}Z(|N_{h}|)\,Z(v)
\\
\nonumber
&=\mnh^{-1}Z(Z(v))+\mnh^{-2}\escpr{N,T}\,(\escpr{B(Z),S}-1)\,Z(v),
\end{align}
where we have used \eqref{eq:vmnh} to compute
$Z(|N_{h}|)$.  To deduce \eqref{eq:divform} it suffices to simplify in
\eqref{eq:toto2} after substituting the information of
\eqref{eq:toto12} and \eqref{eq:toto22}.
\end{proof}

\begin{remark}
If $\Sg$ is a minimal surface then the functional $\mathcal{L}$ in
\eqref{eq:lu} provides a Sturm-Liouville differential operator along
any of the characteristic segments of $\Sg$.
\end{remark}

As a direct consequence of Proposition~\ref{prop:ibp} and
Proposition~\ref{prop:stcond1} we deduce

\begin{corollary}
\label{cor:stcond2}
Let $\Sg$ be a $C^2$ immersed oriented area-stationary surface in
$\hh^1$.  If $\Sg$ is stable under non-singular variations then we
have
\[
-\int_{\Sg}u\,\mathcal{L}(u)\,d\Sg\geq 0,
\]
for any function $u\in C_{0}(\Sg-\Sg_{0})$ which is also $C^2$ in the
$Z$-direction.
\end{corollary}

Finally, with the same technique as in Proposition~\ref{prop:ibp} we
can prove the following lemma.

\begin{lemma}[Integration by parts II]
\label{lem:hordiv}
Let $\Sg\sub\hh^1$ be a $C^2$ immersed surface with unit normal
vector $N$ and singular set $\Sg_{0}$.  Consider two functions $u\in
C_{0}(\Sg-\Sg_{0})$ and $v\in C(\Sg-\Sg_{0})$ which are $C^1$ and
$C^2$ in the $Z$-direction, respectively.  Then we have
\[
\int_{\Sg}\mnh\left\{Z(u)\,Z(v)+u\,Z(Z(v))+
2\,\mnh^{-1}\,\escpr{N,T}\,u\,Z(v)\right\}d\Sg=0.
\]
\end{lemma}

\begin{proof}
Observe that in $\Sg-\Sg_{0}$
\begin{align*}
\divv_{\Sg}(\mnh\,u\,Z(v)\,Z)&=u\,Z(v)\,\big\{Z(\mnh)
+\mnh\,\divv_{\Sg}Z\big\}
\\
&+\mnh\,Z(u)\,Z(v)+\mnh\,u\,Z(Z(v)),
\end{align*}
and that the function in the left-hand side has vanishing integral by
Lemma~\ref{lem:div}.  On the other hand, \eqref{eq:vmnh} gives us
\begin{equation}
\label{eq:znh2}
Z(|N_{h}|)=\escpr{N,T}-\escpr{N,T}\,\escpr{B(Z),S},
\end{equation}
which together with \eqref{eq:toto12} implies
$Z(\mnh)+\mnh \divv_{\Sg}Z=2\,\escpr{N,T}$.  The result follows.
\end{proof}

\begin{remark}
Some other integration by parts formulas in $\hh^1$ can be found in
\cite[Sect.~10]{dgn}.
\end{remark}

\section{Complete stable surfaces with empty singular set}
\label{sec:nosing}
\setcounter{equation}{0}

In this section we provide the classification of $C^2$ complete stable
surfaces in $\hh^1$ with empty singular set.  Recall that if
$\Sg_{0}=\emptyset$ then $\Sg$ is area-stationary if and only if $\Sg$ is minimal by
Proposition~\ref{prop:varprop} (iv).  We say that an immersed surface
$\Sg$ in $\hh^1$ is \emph{complete} if it is complete in the
Riemannian manifold $(\hh^1,g)$.  For a $C^2$ complete area-stationary
surface $\Sg$ with $\Sg_{0}=\emptyset$ the characteristic curves are
straight lines by Proposition~\ref{prop:varprop} (i).  In particular
$\Sg$ cannot be compact.  Some classification results for
area-stationary surfaces with empty singular set can be found in
\cite[Thm.~5.4]{rr1}, \cite{chenghwang} and \cite[Prop.~6.16]{rr2}.
Note also that for such surfaces to be stable is equivalent to be
stable under non-singular variations.

In Euclidean three-space the description of complete stable
area-stationary surfaces can be obtained by means of a logarithmic
cut-off of the function $u=1$ associated to the variation by level
surfaces of the distance function, see \cite{dcp}.  In $\hh^1$ the
vector field induced by the family of equidistants for the
Carnot-Carath\'eodory distance
$d_{cc}$ to a $C^2$ surface with empty singular set coincides, up to a sign,
with the horizontal Gauss map $\nuh$, see \cite[Thm.~1.1 and
1.2]{arcfer}.  This leads us to use the stability condition in
Proposition~\ref{prop:stcond1} with a test function of the form
$f=u\mnh$, where $f$ is continuous with compact support on the surface
and $C^1$ in the direction of the characteristic field $Z$.  We first
compute the index form for these type of functions.

\begin{lemma}
\label{lem:fmnh}
Let $\Sg\sub\hh^1$ be a $C^2$ immersed minimal surface in $\hh^1$ with
unit normal vector $N$ and singular set $\Sg_{0}$.  Then, for any
function $f\in C_{0}(\Sg-\Sg_{0})$ which is also $C^1$ in the
$Z$-direction, we have
\begin{equation}
\label{eq:indexform3}
\mathcal{I}(f\mnh,f\mnh)=\int_{\Sg}\mnh\left\{Z(f)^2-
\mathcal{L}(|N_{h}|)\,f^2\right\}d\Sg,
\end{equation}
where $\mathcal{I}$ is the index form in \eqref{eq:indexform},
and $\mathcal{L}$ is the differential operator in
\eqref{eq:lu}.
\end{lemma}

\begin{proof}
Along this proof we shall denote $w=f\mnh$ and
$q=|B(Z)+S|^2-4|N_{h}|^2$.  Note that $w$ is $C^1$ in the
$Z$-direction and $Z(w)=fZ(\mnh)+\mnh Z(f)$.  If we introduce $w$ in
the index form we obtain
\begin{align}
\label{eq:papa1}
\mathcal{I}(w,w)=\int_{\Sg}\left\{|N_{h}|\,Z(f)^2+|N_{h}|^{-1}
Z(|N_{h}|)^2\,f^2+Z(f^2)\,Z(|N_{h}|)-|N_{h}|\,q\,f^2\right\}d\Sg.
\end{align}
On the other hand, we know from Lemma~\ref{lem:zbzs} (i) that
$|N_{h}|$ is $C^\infty$ in the $Z$-direction.  Therefore, we can apply
Proposition~\ref{prop:ibp} with $u=f^2|N_{h}|$ and $v=|N_{h}|$, so
that we get
\begin{align*}
-\int_{\Sg}|N_{h}|\,\mathcal{L}(|N_{h}|)\,f^2\,d\Sg&=
 \int_{\Sg}\left\{|N_{h}|^{-1}
Z(|N_{h}|)^2\,f^2+Z(f^2)\,Z(|N_{h}|)-|N_{h}|\,q\,f^2\right\}d\Sg
\\
&=\mathcal{I}(w,w)-\int_{\Sg}\mnh\,Z(f)^2\,d\Sg,
\end{align*}
where in the second equality we have used \eqref{eq:papa1}.  This
proves the claim.
\end{proof}

\begin{remark}
Some other versions of \eqref{eq:indexform3} for variations of a $C^2$
surface $\Sg$ with associated vector field $f\nuh$, $f\in
C_{0}^2(\Sg-\Sg_{0})$, can be found in \cite[Lem.~3.9]{dgn3} and
\cite[Thm.~3.4]{dgnp}.  See also \cite[Sect.~3.2]{bscv} and
\cite[Thm.~3.5]{mscv} for the case of an intrinsic graph associated to
a function with less regularity than $C^2$.
\end{remark}

In the next lemma we particularize \eqref{eq:indexform3} for
$f=uv^{-1}$.  This type of test functions will be used
to prove Theorem~\ref{th:nosing}.

\begin{lemma}
Let $\Sg\sub\hh^1$ be a $C^2$ immersed minimal surface in $\hh^1$ with
unit normal vector $N$ and singular set $\Sg_{0}$.  Consider two
functions $u\in C_{0}(\Sg-\Sg_{0})$ and $v\in C(\Sg-\Sg_{0})$ which
are $C^1$ and $C^2$ in the $Z$-direction, respectively.  If $v$ never
vanishes, then
\begin{align}
\label{eq:indexform4b}
\mathcal{I}(uv^{-1}\mnh,&uv^{-1}\mnh)=\int_{\Sg}\mnh\,
v^{-2}Z(u)^2\,d\Sg
\\
\nonumber
&+\int_{\Sg}\mnh\,u^2\bigg\{Z(v^{-1})^2-\frac{1}{2}\,Z(Z(v^{-2}))
-\mnh^{-1}\,\escpr{N,T}\,Z(v^{-2})
\,\bigg\}\,d\Sg
\\
\nonumber
&-\int_{\Sg}\mnh\,\mathcal{L}(\mnh)\,(uv^{-1})^2\,d\Sg,
\end{align}
where $\mathcal{I}$ is the index form in \eqref{eq:indexform},
and $\mathcal{L}$ is the differential operator in
\eqref{eq:lu}.
\end{lemma}

\begin{proof}
From \eqref{eq:indexform3} we only have to compute
\[
\int_{\Sg}\mnh\,Z(uv^{-1})^2\,d\Sg.
\]
Since
\[
Z(uv^{-1})^2=v^{-2}Z(u)^2+u^2Z(v^{-1})^2+
\frac{1}{2}\,Z(u^2)\,Z(v^{-2}),
\]
and Lemma~\ref{lem:hordiv} implies
\[
\int_{\Sg}\frac{1}{2}\,\mnh\,Z(u^2)\,Z(v^{-2})\,d\Sg=
-\int_{\Sg}\mnh\,u^2\bigg\{\frac{1}{2}\,Z(Z(v^{-2}))
+\mnh^{-1}\,\escpr{N,T}\,Z(v^{-2}) \bigg\}\,d\Sg,
\]
we see that \eqref{eq:indexform4b} holds.
\end{proof}

The previous lemmas suggest that, for a function $u=f\mnh$, the
stability condition in Proposition~\ref{prop:stcond1} is more
restrictive if $\mathcal{L}(\mnh)>0$.  Thus it is interesting to
compute $\mathcal{L}(\mnh)$ and to study its sign.

\begin{lemma}
\label{lem:modnh}
Let $\Sg$ be a $C^2$ immersed minimal surface in $\hh^1$ with unit
normal vector $N$.  Consider the basis $\{Z,S\}$ defined in
\eqref{eq:nuh} and \eqref{eq:ese}.  Let $B$ be the Riemannian
shape operator of $\Sg$.  Then, in the regular set $\Sg-\Sg_{0}$,
we have
\begin{equation}
\label{eq:lnh}
\mathcal{L}(|N_{h}|)=4\big(|N_{h}|^{-2}\,\escpr{B(Z),S}-1\big),
\end{equation}
where $\mathcal{L}$ is the second order operator in \eqref{eq:lu}.
\end{lemma}

\begin{proof}
From Lemma~\ref{lem:zbzs} (i) we know that $\mnh$ is $C^\infty$ in the
$Z$-direction.  We must compute $Z(|N_{h}|)$ and $Z(Z(|N_{h}|))$.  By
\eqref{eq:znh2} we have
\[
Z(|N_{h}|)=\escpr{N,T}-\escpr{N,T}\,\escpr{B(Z),S},
\]
and so
\[
Z(Z(|N_{h}|))=Z(\escpr{N,T})-Z(\escpr{N,T})\,\escpr{B(Z),S}
-\escpr{N,T}\,Z(\escpr{B(Z),S}).
\]
Now we use \eqref{eq:znt} and \eqref{eq:zbzs}, so that we get
\begin{align}
\label{eq:zznh}
Z(Z(|N_{h}|))&=-5|N_{h}|+4|N_{h}|^3+2|N_{h}|^{-1}\,\escpr{B(Z),S}
\\
\nonumber
&+2|N_{h}|^{-1}\,\escpr{B(Z),S}^2-3|N_{h}|\,\escpr{B(Z),S}^2.
\end{align}
By substituting \eqref{eq:znh2} and \eqref{eq:zznh} into
\eqref{eq:lu}, we obtain
\begin{align*}
\mathcal{L}(|N_{h}|)&=-5-\escpr{B(Z),S}^2
+4\mnh^{-2}\,\escpr{B(Z),S}-2\escpr{B(Z),S)}+|B(Z)+S|^2
\\
&=4(|N_{h}|^{-2}\,\escpr{B(Z),S}-1),
\end{align*}
where in the second equality we have applied that $\Sg$ is minimal, and
so $B(Z)=\escpr{B(Z),S}S$ by \eqref{eq:mc2}.  This proves
\eqref{eq:lnh}.
\end{proof}

In the next result we show some properties of the Jacobi field
associated to the family of characteristic segments of a minimal
surface in $\hh^1$.  This will allows us to study the sign of
$\mathcal{L}(\mnh)$ and to construct suitable test functions to
introduce in \eqref{eq:indexform4b} when $\Sg$ is a complete minimal
surface with empty singular set.

\begin{lemma}
\label{lem:jacobifield}
Let $\Sg\sub\hh^1$ be a $C^2$ immersed minimal surface with unit
normal $N$ and singular set $\Sg_{0}$.  Consider an integral curve
$\Ga:I\to\Sg-\Sg_{0}$ of the vector field $S$ in \eqref{eq:ese}.  We define
the map $F:I\times I'\to\Sg-\Sg_{0}$ by
$F(\eps,s):=\Ga(\eps)+s\,Z_{\Ga(\eps)}
$.  Let $V_{\eps}(s):=(\ptl F/\ptl\eps)(\eps,s)$.  Then $V_{\eps}$ is
a $C^\infty$ Jacobi vector field along $\ga_{\eps}(s):=F(\eps,s)$.
Moreover, we have
\begin{itemize}
\item[(i)] The vertical component of $V_{\eps}$ is given by
$\escpr{V_{\eps},T}(s)=a_{\eps}s^2+b_{\eps}s+c_{\eps}$, with
\[
b_{\eps}^2-4a_{\eps}c_{\eps}=-\mnh^2(\Ga(\eps))\,
\mathcal{L}(\mnh)(\Ga(\eps)).
\]
\item[(ii)] $V_{\eps}$ is always orthogonal to $\ga_{\eps}$ and
never vanishes along $\ga_{\eps}$.
\item[(iii)] The function
$v_{\eps}(s):=|\escpr{V_{\eps},T}(s)|^{1/2}$ satisfies
\[
Z(v_{\eps}^{-1})^2-\frac{1}{2}\,Z(Z(v_{\eps}^{-2}))-
\mnh^{-1}\escpr{N,T}\,Z(v_{\eps}^{-2})
=\frac{1}{4\,|V_{\eps}|\,\mnh}\,\mathcal{L}(\mnh),
\]
along any segment $\ga_{\eps}(s)$ where $\escpr{V_{\eps},T}(s)$ never
vanishes.
\end{itemize}
\end{lemma}

\begin{proof}
To simplify the notation we will avoid the subscript $\eps$ along the
proof.  We will use primes for both the derivative of functions
depending on $s$ and the covariant derivative along $\ga(s)$.  By
Proposition~\ref{prop:varprop} (i) the curve $\ga$ is a characteristic
curve of $\Sg$.  It follows from \eqref{eq:horgeo} and
Lemma~\ref{lem:jacobi} that $V$ is a $C^\infty$ Jacobi field along
$\ga$ with $[\dot{\ga},V]=0$. Note that
\begin{equation}
\label{eq:prima1}
\escpr{V,T}'=\escpr{V',T}+\escpr{V,T'}=-2\escpr{V,\nuh},
\end{equation}
since $T'=J(Z)=-\nuh$, and
\begin{equation}
\label{eq:vprimat}
\escpr{V',T}=\escpr{D_{Z}V,T}=\escpr{D_{V}Z,T}=
-\escpr{Z,J(V)}=\escpr{J(Z),V}=
-\escpr{V,\nuh}.
\end{equation}
If we derive again in \eqref{eq:prima1} then we obtain
\begin{equation}
\label{eq:prima21}
\escpr{V,T}''=-2\escpr{V',\nuh}-2\escpr{V,\nuh'}=-2\big(\escpr{V',\nuh}
+\escpr{V,T}\big),
\end{equation}
since $\nuh'=D_{Z}\nuh=T$ by \eqref{eq:dznuh} and the fact that $\Sg$
is minimal. Hence
\begin{equation}
\label{eq:prima3}
(-1/2)\escpr{V,T}'''=\big(\escpr{V',\nuh}+\escpr{V,T}\big)'
=\escpr{V'',\nuh}+2\escpr{V',T}+\escpr{V,T'}=0,
\end{equation}
where we have used the Jacobi equation \eqref{eq:jaceq2}, equality
\eqref{eq:vprimat}, and that $T'=-\nuh$.  To simplify
\eqref{eq:prima21} we compute $\escpr{V',\nuh}$.  By \eqref{eq:dvnuh}
and the fact that $V$ is tangent to $\Sg$, we deduce
\begin{align*}
\escpr{V',\nuh}=\escpr{D_{V}Z,\nuh}=-\escpr{Z,D_{V}\nuh}
&=-\mnh^{-1}\,\big(\escpr{D_{V}N,Z}-\escpr{N,T}\,\escpr{J(V),Z}\big)
\\
&=\mnh^{-1}\,\big(\escpr{B(Z),V}+\escpr{N,T}\,\escpr{V,\nuh}\big),
\end{align*}
and so, after substituting into \eqref{eq:prima21}, we get
\begin{equation}
\label{eq:prima22}
\escpr{V,T}''=-2\mnh^{-1}\,\big(\escpr{B(Z),V}+\escpr{N,T}\,\escpr{V,\nuh}
+\mnh\,\escpr{V,T}\big).
\end{equation}

From \eqref{eq:prima3} we conclude that $\escpr{V,T}(s)$ is a
polynomial of degree at most two.  Write
\begin{equation}
\label{eq:poli}
\escpr{V,T}(s)=as^2+bs+c.
\end{equation}
Denote $p=\Ga(\eps)$.  As $V(0)=S_{p}$, it is easy to check by
\eqref{eq:prima1} and \eqref{eq:prima22}, that
\begin{align*}
\nonumber
c&=\escpr{V,T}(0)=-\mnh(p),
\\
\nonumber
b&=\escpr{V,T}'(0)=-2\,\escpr{V,\nuh}(p)=-2\,\escpr{N,T}(p),
\\
\nonumber
a&=(1/2)\,\escpr{V,T}''(0)
=-\mnh^{-1}\big(\escpr{B(Z),S}+\escpr{N,T}^2-\mnh^2\big)(p).
\end{align*}
In particular, it follows from \eqref{eq:lnh}  that
\[
b^2-4ac=-4\,\big(\escpr{B(Z),S}-\mnh^2\big)(p)=
-\mnh^2(p)\,\mathcal{L}(\mnh)(p),
\]
which proves assertion (i) in the statement.

To prove assertion (ii), observe that
\[
\escpr{V,\dot{\ga}}'=\escpr{V',\dot{\ga}}+\escpr{V,\dot{\ga}'}
=\escpr{D_{V}Z,Z}+\escpr{V,D_{Z}Z}=0,
\]
by \eqref{eq:dzz}.  This implies that $\escpr{V,\dot{\ga}}=0$ along
$\ga$ since $V(0)=S_{p}$.  Hence there is a $C^1$ function
$f:I'\to\rr$ such that $V=fS$ along $\ga$.  Clearly $|f|=|V|$, and so
$\escpr{V,T}=\pm |V|\,\mnh$.  By \eqref{eq:poli} the vector $V$
vanishes at most two times along $\ga$.  Suppose that $s_{0}\in I'$ is
the first positive value where $V(s_{0})=0$.  Note that the sign of
$f/|V|$ is constant along a small interval $(s_{0}-\delta,s_{0})$.  By
\eqref{eq:prima1} and \eqref{eq:prima22} we get
$\escpr{V,T}(s_{0})=\escpr{V,T}'(s_{0}) =\escpr{V,T}''(s_{0})=0$.  By
using L'H\^opital's rule twice, we deduce
\begin{align*}
\pm\mnh(\ga(s_{0}))=\lim_{s\uparrow
s_{0}}\,\frac{\escpr{V,T}}{|V|}(s)=\lim_{s\uparrow
s_{0}}\,\frac{|V|\,\escpr{V,T}'}{\escpr{V,V'}}(s)=
\lim_{s\uparrow
s_{0}}\,\frac{|V|'\,\escpr{V,T}'+|V|\,\escpr{V,T}''}
{|V'|^2+\escpr{V,V''}}(s).
\end{align*}
The numerator tends to zero since $|V|'=\escpr{V/|V|,V'}\leq
|V'|\leq M$ on $(s_{0}-\delta,s_{0})$. The denominator
goes to $|V'(s_{0})|^2$, which is positive; otherwise, the
Jacobi field $V$ would be identically zero along $\ga$. It follows
that $\mnh(\ga(s_{0}))=0$, a contradiction since $\ga(s_{0})\in\Sg-\Sg_{0}$.

To prove (iii) let us suppose that $\escpr{V,T}$ never vanishes along
$\ga$.  Then it is clear that
$v=|\escpr{V,T}|^{1/2}=(-\escpr{V,T})^{1/2}$ since $V(0)=S_{p}$.  In
particular, we get $f=|V|>0$ along $\ga$.  Now we derive
$v=(-\escpr{V,T})^{1/2}=(f\mnh)^{1/2}$ with respect to $s$.  By taking
into account \eqref{eq:prima1} and \eqref{eq:prima22}, we obtain
\begin{align*}
Z(v^{-1})&=\frac{1}{2}\,(-\escpr{V,T}^{-3/2})\,\escpr{V,T}'
=-v^{-3}\,\escpr{V,\nuh}=\frac{-\escpr{N,T}}{f^{1/2}\,\mnh^{3/2}},
\\
Z(v^{-2})&=\escpr{V,T}^{-2}\,\escpr{V,T}'
=-2v^{-4}\,\escpr{V,\nuh}=\frac{-2\escpr{N,T}}{f\mnh^2},
\\
Z(Z(v^{-2}))&=\big(\escpr{V,T}^{-2}\escpr{V,T}'\big)'=
-2\escpr{V,T}^{-3}(\escpr{V,T}')^2+\escpr{V,T}^{-2}\escpr{V,T}''
\\
&=\frac{8\escpr{N,T}^2}{f\mnh^3}-\frac{2}{f^2\mnh^3}
\,\big(\escpr{B(Z),V}+\escpr{N,T}\,\escpr{V,\nuh}+\mnh\,\escpr{V,T}\big)
\\
&=\frac{8\escpr{N,T}^2}{f\mnh^3}-\frac{2}{f\mnh^3}
\,\big(\escpr{B(Z),S}+\escpr{N,T}^2 -\mnh^2\big).
\end{align*}
 After simplifying, we conclude by \eqref{eq:lnh} that
\begin{align*}
Z(v^{-1})^2-\frac{1}{2}\,Z(Z(v^{-2}))&-\mnh^{-1}\escpr{N,T}\,Z(v^{-2})
\\
&=\frac{1}{f\mnh}\,\big(\mnh^{-2}\escpr{B(Z),S}-1\big)=
\frac{1}{4f\mnh}\,\mathcal{L}(\mnh),
\end{align*}
which proves the claim.
\end{proof}

\begin{proposition}
\label{lem:lnh>0}
Let $\Sg$ be a $C^2$ complete, oriented, area-stationary surface
immersed in $\hh^1$ with empty singular set.  Then the operator
$\mathcal{L}$ defined in \eqref{eq:lu} satisfies
$\mathcal{L}(\mnh)\geq 0$ on $\Sg$.  Moreover,
$\mathcal{L}(\mnh)(p)=0$ for a point $p\in\Sg$ if and only if
$\escpr{N,T}=0$ and $\escpr{B(Z),S}=1$ along the characteristic line
of $\Sg$ passing through $p$.  As a consequence,
$\mathcal{L}(\mnh)\equiv 0$ on $\Sg$ if and only if any connected
component of $\Sg$ is a Euclidean vertical plane.
\end{proposition}

\begin{proof}
Take a point $p\in\Sg$.  Let $\Ga:I\to\Sg$ be the integral curve
through $p$ of the vector field $S$ in \eqref{eq:ese}.  We define the
map $F:I\times\rr\to\hh^1$ by $F(\eps,s):=\Ga(\eps)+s\,Z_{\Ga(\eps)}$.
By the completeness of $\Sg$ and Proposition~\ref{prop:varprop} (i),
any $\ga_{\eps}(s):=F(\eps,s)$ is a characteristic curve of $\Sg$.  In
particular, $F(I\times\rr)\subseteq\Sg$.

Let $V(s):=(\ptl F/\ptl\eps)(0,s)$.  By using
Lemma~\ref{lem:jacobifield} we deduce that, along the complete line
$\ga(s):=\ga_{0}(s)$, the vectors $V(s)$ and $\dot{\ga}(s)$ generate
the tangent plane to $\Sg$ at $\ga(s)$.  Since $\Sg$ has empty
singular set, it follows that the function $\escpr{V,T}(s)=as^2+bs+c$
never vanishes along $\ga(s)$.  In case $a=0$ we must have $b=0$
(otherwise we would find a root of $as^2+bs+c$).  In case $a\neq 0$ we
must have $b^2-4ac<0$.  Anyway, we get $b^2-4ac\leq 0$ and so
$\mathcal{L}(\mnh)(p)\geq 0$ by Lemma~\ref{lem:jacobifield} (i).

Observe that $\mathcal{L}(\mnh)(p)=0$ if and only if $a=b=0$.  This is
equivalent to that $\escpr{V,T}$ is constant along $\ga$.  It follows
from \eqref{eq:prima1} and \eqref{eq:prima22} that $\escpr{N,T}=0$ and
$\escpr{B(Z),S}=1$ along $\ga$.  Conversely, if $\escpr{N,T}=0$ and
$\escpr{B(Z),S}=1$ along $\ga$ then \eqref{eq:lnh} implies that
$\mathcal{L}(\mnh)=0$ along $\ga$.

Finally, if $\mathcal{L}(\mnh)\equiv 0$ on $\Sg$ then
$\escpr{N,T}\equiv 0$ on $\Sg$.  By \cite[Prop.~6.16]{rr2} we conclude
that any connected component of $\Sg$ must be a Euclidean vertical
plane.  Conversely, it is not difficult to see that
$\mathcal{L}(\mnh)\equiv 0$ holds for any Euclidean vertical plane.
\end{proof}

Now we are ready to prove the main result of this section.

\begin{theorem}
\label{th:nosing}
Let $\Sg$ be a $C^2$ complete, oriented, connected, area-stationary
surface immersed in $\hh^1$ with empty singular set.  If $\Sg$ is not
a Euclidean vertical plane then $\Sg$ is unstable.
\end{theorem}

\begin{proof}
Let $N$ be the unit normal vector to $\Sg$.  We can find $p\in\Sg$
such that $\escpr{N,T}(p)\neq 0$.  Otherwise $\Sg$ would be a
Euclidean vertical plane by \cite[Prop.~6.16]{rr2}.  By using
Proposition~\ref{prop:varprop} (i) and the completeness of $\Sg$, we
can parameterize $\Sg$, around the characteristic line containing $p$,
by the map $F:I\times\rr\to\Sg$ given by
$F(\eps,s)=\Ga(\eps)+s\,Z_{\Ga(\eps)}$, where $\Ga(\eps)$ is a piece
of the integral curve through $p$ of the vector field $S$ in
\eqref{eq:ese}.  Let
$\ga_{\eps}(s):=F(\eps,s)$.  By Lemma~\ref{lem:jacobifield} we know
that $V_{\eps}(s):=(\ptl F/\ptl\eps)(\eps,s)$ is a non-vanishing
Jacobi field orthogonal to $\ga_{\eps}(s)$.  Moreover, the function
$\escpr{V_{\eps}(s),T}$ is strictly negative since $\Sg$ has empty
singular set and $V_{\eps}(0)=S_{\Ga(\eps)}$.  We consider the
function
$v(\eps,s):=|\escpr{V_{\eps}(s),T}|^{1/2}=(\mnh\,|V_{\eps}(s)|)^{1/2}$,
which is continuous and $C^\infty$ along any $\ga_{\eps}(s)$.

Now we use the coarea formula to compute the index form
\eqref{eq:indexform4b} in terms of the coordinates $(\eps,s)$.  The
Riemannian area element can be expressed as
\[
d\Sg=|V_{\eps}|\,d\eps\,ds.
\]
Hence by using the definition of $v$ together with
Lemma~\ref{lem:jacobifield} (iii), equation \eqref{eq:indexform4b} reads
\begin{equation}
\label{eq:indexfinal}
\mathcal{I}(uv^{-1}\mnh,uv^{-1}\mnh)=\int_{I\times\rr}
\left(\frac{\ptl u}{\ptl s}\right)^2d\eps\,ds
-\frac{3}{4}\int_{I\times\rr}\mathcal{L}(\mnh)\,u^2\,d\eps\,ds,
\end{equation}
for any $u\in C_{0}(I\times\rr)$ which is also $C^1$ with
respect to $s$.

Take a non-negative $C^\infty$ function $\phi:I\to\rr$ with
$\phi(0)>0$ and compact support contained inside a bounded interval
$I'\subseteq I$.  Denote $\ell:=\text{length}(I')$.  Let $M$ be a
positive constant so that $|\phi'(\eps)|\leq M$,  $\eps\in I$.
For any $k\in\mathbb{N}$ we define the function
\[
u_{k}(\eps,s):=\phi(\eps)\,\phi(s/k).
\]
It is clear that $u_{k}\in C_{0}(I'\times kI')$, and that $u_{k}$ is
$C^\infty$ with respect to $s$.  By Fubini's theorem
\[
\int_{I\times\rr}\left(\frac{\ptl u_{k}}{\ptl s}\right)^2d\eps\,ds=
\frac{1}{k^2}\left(\int_{I'}\phi(\eps)^2\,d\eps\right)
\,\left(\int_{kI'}\phi'(s/k)^2\,ds\right)\leq\frac{\ell M^2}{k}
\int_{I'}\phi(\eps)^2\,d\eps,
\]
which goes to $0$ when $k\to\infty$.  Note also that
$\{u_{k}\}_{k\in\mathbb{N}}$ pointwise converges when $k\to\infty$ to
$u(\eps,s)=\phi(0)\,\phi(\eps)$.  By Proposition~\ref{lem:lnh>0} we have
$\mathcal{L}(\mnh)\geq 0$ on $\Sg$.  Thus we can apply Fatou's lemma
to obtain
\[
\liminf_{k\to\infty}\int_{I\times\rr}\mathcal{L}(\mnh)\,u^2_{k}\,d\eps\,ds
\geq\int_{I\times\rr}\mathcal{L}(\mnh)\,u^2\,d\eps\,ds.
\]
We conclude from \eqref{eq:indexfinal} that
\begin{align*}
\limsup_{k\to\infty}\,\mathcal{I}(u_{k}v^{-1}\mnh,u_{k}v^{-1}\mnh)&=
-\frac{3}{4}\,\liminf_{k\to\infty}\int_{I\times\rr}\mathcal{L}(\mnh)\,u^2_{k}\,
d\eps\,ds
\\
\nonumber
&\leq-\frac{3}{4}\int_{I\times\rr}\mathcal{L}(\mnh)\,
u^2\,d\eps\,ds,
\end{align*}
which is strictly negative by Proposition~\ref{lem:lnh>0} since
$\escpr{N,T}\neq 0$ inside an open neighborhood around $p$.  Hence
$\Sg$ is unstable.
\end{proof}

\begin{corollary}
\label{cor:nosing}
Let $\Sg$ be a $C^2$ complete, oriented, connected, area-stationary
surface immersed in $\hh^1$ with empty singular set.  Then $\Sg$ is
stable if and only if $\Sg$ is a Euclidean vertical plane.
\end{corollary}

\begin{proof}
The necessary condition follows from Theorem~\ref{th:nosing}.
Conversely, suppose that $\Sg$ is a vertical plane.  We can prove that
$\Sg$ is an area-minimizing surface in $\hh^1$ by using a
calibration argument similar to the one in \cite[Thm.~5.3]{rr2}, see
also \cite[Ex.~2.2]{bscv}. In particular, $\Sg$ is stable.
\end{proof}

\begin{remark}
Previous results related to Corollary~\ref{cor:nosing} were obtained
in \cite{bscv} and \cite{dgnp}.  Precisely, in \cite[Thm.~5.1]{bscv}
it is proved that the Euclidean vertical planes are the only complete
stable intrinsic graphs in $\hh^1$ associated to a $C^2$ function.  In
\cite[Thm.~1.8]{dgnp} vertical planes are characterized as the unique
complete stable $C^2$ Euclidean graphs with empty singular set.  As we
pointed out in the introduction of the paper,
Corollary~\ref{cor:nosing} does not follow from the aforementioned
results.  For example, they do not apply for the family of
sub-Riemannian catenoids $t^2=\la^2\,(x^2+y^2-\la^2)$, $\la\neq 0$.
\end{remark}

\section{Complete stable surfaces with non-empty singular set}
\label{sec:singcurves}
\setcounter{equation}{0}

In this section we give the classification of $C^2$ complete stable
surfaces in $\hh^1$ with non-empty singular set.  By
Proposition~\ref{prop:varprop} the singular set of a $C^2$
area-stationary surface consists of isolated points and curves of
class $C^1$.  Moreover, the characteristic curves in the regular set
meet the singular curves orthogonally.  By using these facts we
were able to obtain the following result in \cite[Thm.~6.15]{rr2}.

\begin{proposition}
\label{prop:basic2}
Let $\Sg$ be a $C^2$ complete, oriented, connected, area-stationary
surface immersed in $\hh^1$ with singular set $\Sg_{0}$.
\begin{itemize}
\item[(i)] If $\Sg_{0}$ contains an isolated point then $\Sg$
coincides with a Euclidean non-vertical plane.  \item[(ii)] If
$\Sg_{0}$ contains a singular curve then $\Sg$ is either congruent to
the hyperbolic paraboloid $t=xy$ or to one of the helicoidal surfaces
$\h_{R}$ defined below.
\end{itemize}
\end{proposition}

In \cite[Ex.~6.14]{rr2} we described the helicoid $\h_{R}$ as the
union of all the horizontal straight lines orthogonal to the
sub-Riemannian geodesic in $\hh^1$ obtained by the horizontal lift of the
circle in the $xy$-plane of radius $1/R$ centered at the origin.  We
can parameterize $\h_{R}$ by means of the $C^\infty$ diffeomorphism
$F:\rr^2\to\h_{R}$ defined by
\begin{equation}
\label{eq:efeh}
F(\eps,s)=(s\,\sin(R\eps),s\,\cos(R\eps),\eps/R).
\end{equation}
The singular set of $\h_{R}$ consists of the helices $s=\pm 1/R$.
Note that the family $\{\h_{R}\}_{R>0}$ is invariant under the
dilations $\delta_{\la}$ defined in \eqref{eq:dilations}.  In fact, it
can be checked from \eqref{eq:efeh} that
$\delta_{\la}(\h_{R})=\h_{{R'}}$ with $R'=e^{-\la}R$.  The surfaces
$\h_{R}$ coincide with the classical left-handed minimal helicoids in
$\rr^3$.  In particular, they are embedded surfaces containing the
vertical axis.  We remark that the classical right-handed minimal
helicoids in $\rr^3$ are complete area-stationary surfaces in $\hh^1$
with empty singular set, and so they are unstable by
Theorem~\ref{th:nosing}.

Proposition~\ref{prop:basic2} indicates us that the study of stable
surfaces in $\hh^1$ with non-empty singular set can be reduced to
three cases: Euclidean non-vertical planes, the hyperboloid $t=xy$ and
the helicoids $\h_{R}$.  In \cite[Thm.~5.3]{rr2} we showed that any
complete $C^2$ area-stationary graph over the $xy$-plane is an area-minimizing
surface.  This gives us the stability of any plane $t=ax+by$ and any
surface congruent to $t=xy$.  So it remains to analyze the stability
of the helicoidal surfaces $\h_{R}$.

We first compute some geometric terms of a helicoid $\h_R$ with
respect to the system of coordinates $(\eps,s)$ in \eqref{eq:efeh}.
Note that
\begin{align}
\nonumber
\frac{\ptl F}{\ptl\eps}&=Rs\,\cos(R\eps)\,X-Rs\,\sin(R\eps)\,Y+
f(s)\,T,
\\
\nonumber
\frac{\ptl F}{\ptl s}&=\sin(R\eps)\,X+\cos(R\eps)\,Y,
\end{align}
where $f:\rr\to\rr$ is defined by
\[
f(s)=\frac{1}{R}-Rs^2.
\]
As a consequence, the Riemannian area element is given by
\begin{equation}
\label{eq:jacfr}
d\Sg=\sqrt{f(s)^2+R^2s^2}\,\,d\eps\,ds.
\end{equation}
On the other hand, the cross product of $\ptl F/\ptl s$ and
$\ptl F/\ptl\eps$ in $(\hh^1,g)$ provides the following unit normal
vector to $\h_R$
\begin{equation}
\label{eq:nhr}
N=\frac{f(s)\,\cos(R\eps)\,X-f(s)\,\sin(R\eps)\,Y-Rs\,T}
{\sqrt{f(s)^2+R^2s^2}},
\end{equation}
and so
\begin{equation}
\label{eq:xopo2}
\mnh=\frac{|f(s)|}{\sqrt{f(s)^2+R^2s^2}},\qquad
\escpr{N,T}=\frac{-Rs}{\sqrt{f(s)^2+R^2s^2}}.
\end{equation}
It follows that the straight lines $\ga_{\eps}(s)=F(\eps,s)$,
$s\in\rr$, satisfy
\begin{equation}
\label{eq:sign}
\dot{\ga}_{\eps}(s)=\text{sign}(1/R-|s|)\,Z,\qquad
|s|\neq 1/R.
\end{equation}
By taking into account \eqref{eq:znt} and \eqref{eq:xopo2} we get,
for $|s|\neq 1/R$, that
\begin{equation}
\label{eq:bzs2}
\escpr{B(Z),S}=1+\mnh^{-1}\,Z(\escpr{N,T})
=\frac{2f(s)^2-Rf(s)}{f(s)^2+R^2s^2}-1,
\end{equation}
which in particular implies
\begin{equation}
\label{eq:qhr}
|B(Z)+S|^2-4\mnh^2=\frac{(R^2-4)\,f(s)^2}{(f(s)^2+R^2s^2)^2}.
\end{equation}

Now we are ready to deduce from Theorem~\ref{th:2ndvar}
and Proposition~\ref{prop:2ndvert} a stability criterion for helicoidal
surfaces that plays the same role as Proposition~\ref{prop:stcond1}.

\begin{proposition}
\label{prop:2ndvarsing}
Let $\Sg$ be the helicoid $\h_R$. If $\Sg$ is stable then, for any function
$u\in C^2_0(\Sg)$ such that $Z(u)=0$ inside a small tubular neighborhood
of $\Sg_0$, we have  $\mathcal{Q}(u)\geq 0$, where
\begin{align}
\label{eq:Qu}
\mathcal{Q}(u):=&\int_{\Sg}\mnh^{-1}\left\{Z(u)^2-\big(|B(Z)+S|^2
-4\mnh^2\big)\,u^2\right\}d\Sg
\\
\nonumber
&-4\int_{\Sg_{0}}u^2\,d\Sg_{0}+\int_{\Sg_0}S(u)^2\,d\Sg_0.
\end{align}
Here $\{Z,S\}$ is the orthonormal basis in \eqref{eq:nuh} and \eqref{eq:ese},
$B$ is the Riemannian shape operator of $\Sg$, and $d\Sg_{0}$ is the
Riemannian length measure on $\Sg_0$.
\end{proposition}

\begin{proof}
We suppose that the unit normal $N$ to $\Sg$ is the one in \eqref{eq:nhr}. For simplicity we denote $q=|B(Z)+S|^2-4\mnh^2$. By \eqref{eq:qhr}, \eqref{eq:xopo2} and \eqref{eq:jacfr} it follows that $\mnh^{-1}qu^2\in L^1(\Sg)$ provided $u\in C_0(\Sg)$. In particular, $\mathcal{Q}(u)$ is well defined for any $u\in C_0(\Sg)$ which is piecewise $C^1$ in the $Z$-direction, satisfies $\mnh^{-1}Z(u)^2\in L^1(\Sg)$, and whose restriction to $\Sg_0$ is $C^1$.

Let us show, in a first step, the following statement
\begin{align}
\label{eq:inehr}
\mathcal{Q}(v)\geq 0, \  &\text{for any } v\in C_0^2(\Sg) \text{ such that }
Z(v/\escpr{N,T})=0
\\
\nonumber
& \text{ in a small tubular neighborhood } E \text{ of } \Sg_0.
\end{align}
Note that a function $v$ as above satisfies $Z(v)^2=(Z(\escpr{N,T})^2/\escpr{N,T}^2)\,v^2$ in $E$. It follows from \eqref{eq:xopo2} and \eqref{eq:jacfr} that $\mnh^{-1}Z(v)^2\in L^1(\Sg)$, and so $\mathcal{Q}(v)<\infty$.

Let $\sg_0$ be the radius of $E$ and $K$ the support of $v$. For any $\sg\in (0,\sg_0/2)$ let $E_\sg$ be the tubular neighborhood of $\Sg_0$ of radius $\sg$. We consider functions $h_\sg,g_\sg\in C^\infty_0(\Sg)$ such that $g_\sg=1$ on $K\cap\overline{E}_\sg$, $\text{supp}(g_\sg)\subset E_{2\sg}$ and $h_\sg+g_\sg=1$ on $K$. We define the $C^2$ vector field
\[
U_\sg:=(h_\sg\,v)\,N+g_\sg\,\frac{v}{\escpr{N,T}}\,T,
\]
whose support is contained in $K$. Note that $\escpr{U_\sg,N}=v$ on $K$. Let $\varphi_r^\sg(p):=\exp_p(r (U_\sg)_p)$ be the variation associated to $U_\sg$ and $A_\sg(r):=A(\varphi_r^\sg(\Sg))$ the corresponding area functional. The variation $\var_r^\sg$ is vertical when restricted to $E_\sg$. Hence we can suppose, by applying Proposition~\ref{prop:2ndvert} to $w=v/\escpr{N,T}$, that the second derivative of $A_{1\sg}(r):=A(\varphi_r^\sg(E_\sg))$ is given by
\[
A_{1\sg}''(0)=\int_{\Sg_0}S(v)^2\,d\Sg_0.
\]
In the previous equality we have used that $\escpr{N,T}=\pm 1$ on $\Sg_0$. On the other hand, the second derivative of $A_{2\sg}(r):=A(\varphi_r^\sg(\Sg-E_\sg))$ can be computed from Theorem~\ref{th:2ndvar}. We obtain the following expression
\begin{align*}
A_{2\sg}''(0)&=\int_{\Sg-E_\sg}\mnh^{-1}\left\{Z(v)^2-qv^2\right\}d\Sg
+\int_{\Sg-E_\sg}\divv_\Sg(\xi Z)\,d\Sg+\int_{\Sg-E_\sg}\divv_\Sg(\mu Z)\,d\Sg,
\end{align*}
where
\begin{align*}
\xi&=\escpr{N,T}\,(1-\escpr{B(Z),S})\,v^2,
\\
\mu&=\mnh^2\left\{\escpr{N,T}\,(1-\escpr{B(Z),S})\,\,
\frac{g^2_\sg\,v^2}{\escpr{N,T}^2}
-2\escpr{B(Z),S}\,\,\frac{h_\sg\,g_\sg\,v^2}{\escpr{N,T}}\right\}.
\end{align*}
If $\Sg$ is stable then $A_\sg''(0)\geq 0$. As $A_\sg(r)=A_{1\sg}(r)+A_{2\sg}(r)$ we deduce, by using the classical Riemannian divergence theorem that, for any $\sg\in (0,\sg_0/2)$, we have the inequality
\begin{equation}
\label{eq:d2asg}
\int_{\Sg-E_\sg}\mnh^{-1}\left\{Z(v)^2-qv^2\right\}d\Sg
-\int_{\ptl E_\sg}(\xi+\mu)\,\escpr{Z,\eta}\,dl+\int_{\Sg_0}S(v)^2\,d\Sg_0\geq 0,
\end{equation}
where $\eta$ is the unit normal to $\ptl E_\sg$ pointing into $\Sg-E_\sg$ and $dl$ denotes the Riemannian length element.

Let us compute the boundary term above. Fix $k\in\{1,2\}$. Let $\Lambda$ be one of the two components of $\ptl E_\sg$ at distance $\sg$ of the singular curve where $\escpr{N,T}=(-1)^{k+1}$. By taking into account \eqref{eq:sign} it follows that $\eta=(-1)^{k+1}Z$ along $\Lambda$. Moreover, the functions $\xi$ and $\mu$ are constant along $\Lambda$. Since $g_\sg=1$ and $h_\sg=0$ on $\Lambda$ we have
\begin{align}
\label{eq:extrahr1}
\int_\Lambda\xi\,\escpr{Z,\eta}\,d\Lambda&=
(-1)^{k+1}\,\escpr{N,T}\,(1-\escpr{B(Z),S})\,
\int_\Lambda v^2\,d\Lambda,
\\
\label{eq:extrahr2}
\int_\Lambda\mu\,\escpr{Z,\eta}\,d\Lambda&=
(-1)^{k+1}\,\mnh^2\,\escpr{N,T}^{-1}\,(1-\escpr{B(Z),S})\,
\int_\Lambda v^2\,d\Lambda.
\end{align}

Now we let $\sg\to 0$ in \eqref{eq:d2asg}. From the dominated convergence
theorem we get
\begin{align*}
\lim_{\sg\to 0}\,\int_{\Sg-E_\sg}\mnh^{-1}\left\{Z(v)^2-qv^2\right\}d\Sg&=
\int_{\Sg}\mnh^{-1}\left\{Z(v)^2-qv^2\right\}d\Sg,
\\
\lim_{\sg\to 0}\,\int_{\ptl E_\sg}v^2\,dl&=2
\int_{\Sg_0}v^2\,d\Sg_0.
\end{align*}
On the other hand, equation \eqref{eq:bzs2} yields $\escpr{B(Z),S}\to -1$ when we approach $\Sg_0$. Moreover, we know that $\mnh\to 0$ and $\escpr{N,T}\to\pm 1$ when $\sg\to 0$. This facts, together with \eqref{eq:extrahr1} and \eqref{eq:extrahr2}  imply that
\[
\lim_{\sg\to 0}\,\int_{\ptl E_\sg}(\xi+\mu)\,\escpr{Z,\eta}\,dl=
4\int_{\Sg_0}v^2\,d\Sg_0.
\]
Hence we obtain $\mathcal{Q}(v)\geq 0$ from \eqref{eq:d2asg}. This proves \eqref{eq:inehr}.

Now we take $u\in C^2_0(\Sg)$ with $Z(u)=0$ inside a small tubular neighborhood $E$ of $\Sg_0$. For any $\sg\in (0,1)$ let $D_\sg$ be the open neighborhood of $\Sg_0$ such that $|\escpr{N,T}|=1-\sg$ on $\ptl D_\sg$. We can find $\sg_0>0$ such that $D_\sg\subset E$ for $\sg\in (0,\sg_0)$. For such values of $\sg$ we define the function $\phi_\sg:\Sg\to [0,1]$ given by
\[
\phi_\sg=
\begin{cases}
|\escpr{N,T}|,\quad\text{in } \overline{D}_\sg,
\\
1-\sg, \qquad\,\,\text{in } \Sg-D_\sg.
\end{cases}
\]
Clearly $\phi_\sg$ is continuous and piecewise $C^1$ in the $Z$-direction. Moreover, the sequence $\{\phi_\sg\}_{\sg\in (0,\sg_0)}$ pointwise converges to $1$ when $\sg\to 0$. By using \eqref{eq:xopo2} and \eqref{eq:jacfr} we can see that $\mnh^{-1}Z(\escpr{N,T})^2$ extends to a continuous function on $\Sg$, and so
\[
\lim_{\sg\to 0}\,\int_\Sg\mnh^{-1}\,Z(\phi_\sg)^2\,d\Sg=0.
\]
By a standard approximation argument we can slightly modify $\phi_\sg$ around $\ptl D_\sg$ in order to construct a sequence of $C^2$ functions $\{\psi_\sg\}_{\sg\in(0,\sg_0)}$ satisfying the same properties. Define $v_\sg:=\psi_\sg u$. This provides a sequence of functions in $C^2_0(\Sg)$ such that $v_\sg=u$ in $\Sg_0$ and $Z(v_\sg/\escpr{N,T})=0$ inside a small tubular neighborhood of $\Sg_0$.  As a consequence of \eqref{eq:inehr} we have $\mathcal{Q}(v_\sg)\geq 0$ for any $\sg\in (0,\sg_0)$. Finally, it is straightforward to check by using the dominated convergence theorem and the Cauchy-Schwartz inequality in $L^2(\Sg)$ that $\{\mathcal{Q}(v_\sg)\}\to\mathcal{Q}(u)$ when $\sg\to 0$. The proposition is proved.
\end{proof}

\begin{remark}
By using Proposition~\ref{prop:2ndvert}, inequality $\mathcal{Q}(u)\geq 0$ can be generalized for any $C^2$ stable solution $\Sg$ of the Plateau problem whose singular curves are of class $C^3$ whenever $\supp(u)$ is contained in the interior of $\Sg$
and $Z(u)=0$ inside a tubular neighborhood of the singular curves.
\end{remark}

Now we are ready to prove the main result of this section.

\begin{theorem}
\label{th:helicoids}
The helicoidal surfaces $\h_{R}$ are all unstable.
\end{theorem}

\begin{proof}
To prove the claim it suffices to show that $\h_{2}$ is unstable.  In
fact, for any $R>0$ we have $\h_{R}=\delta_{\la}(\h_{2})$, where
$\delta_{\la}$ is the dilation defined in \eqref{eq:dilations} with
$\la=\log(2/R)$. By virtue of Lemma~\ref{lem:dilinv} we deduce that
$\h_{R}$ is stable if and only if $\h_{2}$ is stable.

Let $\Sg:=\h_{2}$.  Consider the diffeomorphism $F:\rr^2\to\Sg$ in
\eqref{eq:efeh}.  We denote $\ga_{\eps}(s)=F(\eps,s)$, $s\in\rr$.  The
singular set $\Sg_{0}$ consists of the singular curves
$F(\eps,-1/2)$ and $F(\eps,1/2)$, $\eps\in\rr$.  We
suppose that the normal $N$ to $\Sg$ is the one in
\eqref{eq:nhr}. By  equation \eqref{eq:qhr} we get
\[
|B(Z)+S|^2-4\mnh^2=0,\qquad\text{ on } \Sg-\Sg_{0}.
\]
In particular, the quadratic form $\mathcal{Q}$ in
\eqref{eq:Qu} is given by
\begin{equation}
\label{eq:d2ah2}
\mathcal{Q}(u)=
\int_{\Sg}\mnh^{-1}\,Z(u)^2\,d\Sg-4\int_{\Sg_0}u^2\,d\Sg_0
+\int_{\Sg_0}S(u)^2\,d\Sg_0,
\end{equation}
for any $u\in C_0(\Sg)$ which is piecewise $C^1$ in the
$Z$-direction, satisfies $\mnh^{-1}Z(u)^2\in L^1(\Sg)$, and
whose restriction to $\Sg_0$ is $C^1$.  We apply in \eqref{eq:d2ah2}
the coarea formula.  By using \eqref{eq:xopo2}, \eqref{eq:sign}
and \eqref{eq:jacfr}, we deduce that
\begin{align}
\label{eq:if5}
\mathcal{Q}(u)&=\int_{\rr^2}\frac{f(s)^2+4s^2}{|f(s)|}
\,\left(\frac{\ptl u}{\ptl s}\right)^2d\eps\,ds-4\int_{\rr}u(\eps,-1/2)^2\,d\eps
-4\int_{\rr}u(\eps,1/2)^2\,d\eps
\\
\nonumber
&+\int_\rr\left(\frac{d}{d\eps}\,u(\eps,-1/2)\right)^2d\eps
+\int_\rr\left(\frac{d}{d\eps}\,u(\eps,1/2)\right)^2d\eps.
\end{align}

Let $\phi:\rr\to\rr$ be any $C^\infty$ function with compact support
$[-\eps_0,\eps_0]$. For any $k>1/2$ and $\delta>0$, let
$\phi_{k\delta}:\rr\to [0,1]$ be the symmetric function with respect to
the origin given, for $s\geq 0$, by
\[
\phi_{k\delta}(s)=
\left\{
\begin{array}{ll}
1,& 0\leq s\leq k,\\
\displaystyle\delta^{-1}\,(-s+\delta+k),& k\leq s\leq k+\delta,\\
0,&  s\geq k+\delta.\\
\end{array}
\right.
\]
Now we define the function $u_{k\delta}$ on $\Sg$ whose expression
in coordinates $(\eps,s)$ is
\[
u_{k\delta}(\eps,s)=\phi(\eps)\,\phi_{k\delta}(s).
\]
Clearly $u_{k\delta}$ is a function in $C_{0}(\Sg)$ which is also
$C^\infty$ with respect to $\eps$ and piecewise $C^\infty$ in the
$Z$-direction. Note also that
\begin{equation}
\label{eq:valuescurves}
u_{k\delta}(\eps,-1/2)=u_{k\delta}(\eps,1/2)=\phi(\eps),
\qquad\eps\in\rr.
\end{equation}
Moreover $(\ptl u_{k\delta}/\ptl
s)(\eps,s)=\phi(\eps)\,\phi_{k\delta}'(s)$, which vanishes if $|s|<k$
or $|s|>k+\delta$, and equals $\pm\phi(\eps)/\delta$ if
$k<|s|<k+\delta$. This implies that $Z(u_{k\delta})=0$ inside a
tubular neighborhood of $\Sg_0$. By using Fubini's theorem and
that $|f(s)|^{-1}(f(s)^2+4s^2)$ is symmetric with respect to the
origin, we have
\begin{equation}
\label{eq:mathe}
\int_{\rr^2}\frac{f(s)^2+4s^2}{|f(s)|}
\,\left(\frac{\ptl u_{k\delta}}{\ptl s}\right)^2d\eps\,ds
=\left(\int_{-\eps_0}^{\eps_0}\phi(\eps)^2\,d\eps\right)\,
\left (\frac{2}{\delta^2}\int_{k}^{k+\delta}
\frac{f(s)^2+4s^2}{|f(s)|}\,ds\right).
\end{equation}
The second integral in the right-hand side can be easily
computed. We obtain
\[
2\int_{k}^{k+\delta}
\frac{f(s)^2+4s^2}{|f(s)|}\,ds=\int_{k}^{k+\delta}
\frac{16 s^4+8s^2+1}{4s^2-1}\,ds=\frac{4s^3}{3}+3s
+\log\left(\frac{2s-1}{2s+1}\right)\,\bigg]_{k}^{k+\delta}.
\]
By an elementary analysis we can find a value $k>1/2$ and
$\delta=2k+1$ such that the integral above times $1/\delta^2$ is strictly less
than $8$.  By substituting this information into \eqref{eq:mathe}, and
using \eqref{eq:if5} together with \eqref{eq:valuescurves}, we
conclude for $v:=u_{k\delta}$
\[
\mathcal{Q}(v)<M\,\int_{-\eps_0}^{\eps_0}\phi(\eps)^2\,d\eps
+2\int_{-\eps_0}^{\eps_0}\phi'(\eps)^2\,d\eps,
\]
for some constant $M<0$ which does not depend on the function $\phi$.
If $\eps_0$ is large enough, then we can choose $\phi$
with compact support $[-\eps_0,\eps_0]$ such that the right-hand side
of the previous equation is strictly negative. This can be done since
\[
\inf\left\{\left(\int_\rr\phi'(\eps)^2\,d\eps\right)
\left(\int_\rr\phi(\eps)^2\,d\eps\right)^{-1};\phi\in C_0^\infty(\rr)\right\}=0.
\]

Denote $\bar{\phi}:=\phi_{k\delta}$ for the particular values
of $k$ and $\delta$ found above. We mollify $\bar{\phi}$ in
order to obtain a sequence of functions
$v_\sg(\eps,s)=\phi(\eps)\bar{\phi}_{\sg}(s)$ in $C^\infty_0(\Sg)$
with $v_\sg=v$ on $\Sg_0$ and
\[
\lim_{\sg\to 0}\,\int_\Sg\mnh^{-1}Z(v_\sg)^2\,d\Sg=
\int_\Sg\mnh^{-1}Z(v)^2\,d\Sg.
\]
Hence we have
\[
\lim_{\sg\to 0}\mathcal{Q}(v_\sg)=\mathcal{Q}(v)<0.
\]
By Proposition~\ref{prop:2ndvarsing} we conclude that
$\Sg$ is unstable.
\end{proof}

\begin{remark}
Though the helicoids $\h_{R}$ are unstable, it is possible to
obtain by means of a calibration argument similar to the one used for
the hyperboloid $t=xy$ in \cite[Thm.~5.3]{rr2} that the surface
obtained by removing the vertical axis from $\h_{R}$ is
area-minimizing. On the other hand, the second derivative of the area
in Theorem~\ref{th:2ndvar} indicates us that
any non-singular variation induced by a vector field $U=vN+wT$ such that $v$
and $w$ are $C^1$ functions whose support is contained in the regular set of
$\h_{2}$ satisfies $A''(0)\geq 0$.  This means that $\h_{2}$ is also
\emph{stable under the variations used in} Theorem~\ref{th:2ndvar}. The proof
of Theorem~\ref{th:helicoids} shows that, to get that $\h_2$ is unstable, we need to
consider a function whose support intersects a large piece of $\h_2$ containing the
vertical axis and the singular set.
\end{remark}

\section{Main result}
\label{sec:mainresult}

As a consequence of our previous stability results we can prove the
following.

\begin{theorem}
\label{th:main}
Let $\Sg$ be a $C^2$ complete, oriented, connected, area-stationary
surface immersed in $\hh^1$.  Then $\Sg$ is stable if and only if
$\Sg$ is a Euclidean plane or $\Sg$ is congruent to the
hyperbolic paraboloid $t=xy$.  In particular, $\Sg$ is
area-minimizing.
\end{theorem}

\begin{proof}
If $\Sg$ is stable and the singular set $\Sg_{0}$ is empty then $\Sg$
must be a vertical plane by Theorem~\ref{th:nosing}.  If $\Sg$ is
stable and $\Sg_{0}\neq\emptyset$ then Proposition~\ref{prop:basic2}
and Theorem~\ref{th:helicoids} imply that $\Sg$ coincides with a
non-vertical Euclidean plane, or it is congruent to the hyperbolic
paraboloid $t=xy$.  That Euclidean planes and surfaces congruent to
$t=xy$ are area-minimizing follows from \cite[Ex.~2.2]{bscv} and
\cite[Thm.~5.3]{rr2}.
\end{proof}

\providecommand{\bysame}{\leavevmode\hbox
to3em{\hrulefill}\thinspace}
\providecommand{\MR}{\relax\ifhmode\unskip\space\fi MR }
\providecommand{\MRhref}[2]{%
  \href{http://www.ams.org/mathscinet-getitem?mr=#1}{#2}
} \providecommand{\href}[2]{#2}


\begin{thebibliography}{10}

\bibitem{arcfer}
N.~Arcozzi and F.~Ferrari, \emph{Metric normal and distance
function in the
  {H}eisenberg group}, Math. Z. \textbf{256} (2007), no.~3, 661--684.
  \MR{MR2299576}

\bibitem{balogh}
Z.~M. Balogh, \emph{Size of characteristic sets and functions with
prescribed
  gradient}, J. Reine Angew. Math. \textbf{564} (2003), 63--83. \MR{MR2021034
  (2005d:43007)}

\bibitem{bdce}
J.~L. Barbosa, M.~P. do~Carmo, and J.~Eschenburg, \emph{Stability
of
  hypersurfaces of constant mean curvature in {R}iemannian manifolds}, Math. Z.
  \textbf{197} (1988), no.~1, 123--138. \MR{MR917854 (88m:53109)}

\bibitem{bscv}
V.~Barone~Adesi, F.~Serra~Cassano, and D.~Vittone, \emph{The
{B}ernstein
  problem for intrinsic graphs in {H}eisenberg groups and calibrations}, Calc.
  Var. Partial Differential Equations \textbf{30} (2007), no.~1, 17--49.
  \MR{MR2333095}

\bibitem{andre}
A.~Bella{\"{\i}}che, \emph{The tangent space in sub-{R}iemannian
geometry},
  Sub-Riemannian geometry, Progress in Mathematics, vol. 144, Birkh\"auser,
  Basel, 1996, pp.~1--78. \MR{MR1421822 (98a:53108)}

\bibitem{bernstein}
S.~Bernstein, \emph{{Sur un th{\'e}or{\`e}me de g{\'e}om{\'e}trie
et son
  application aux {\'e}quations aux d{\'e}riv{\'e}es partielles du type
  elliptique}}, Charikov, Comm. Soc. Math. (2) \textbf{15} (1915-1917), 38--45
  (French).

\bibitem{blair}
D.~E. Blair, \emph{Riemannian geometry of contact and symplectic
manifolds},
  Progress in Mathematics, vol. 203, Birkh\"auser Boston Inc., Boston, MA,
  2002. \MR{MR1874240 (2002m:53120)}

\bibitem{survey}
L.~Capogna, D.~Danielli, S.~D. Pauls, and J.~T. Tyson, \emph{An
introduction to
  the {H}eisenberg group and the sub-{R}iemannian isoperimetric problem},
  Progress in Mathematics, vol. 259, Birkh\"auser Verlag, Basel, 2007.
  \MR{MR2312336}

\bibitem{chenghwang}
J.-H. Cheng and J.-F. Hwang, \emph{Properly embedded and immersed
minimal
  surfaces in the {H}eisenberg group}, Bull. Austral. Math. Soc. \textbf{70}
  (2004), no.~3, 507--520. \MR{MR2103983 (2005f:53010)}

\bibitem{chmy}
J.-H. Cheng, J.-F. Hwang, A.~Malchiodi, and P.~Yang, \emph{Minimal
surfaces in
  pseudohermitian geometry}, Ann. Sc. Norm. Super. Pisa Cl. Sci. (5) \textbf{4}
  (2005), no.~1, 129--177. \MR{MR2165405 (2006f:53008)}

\bibitem{chy}
J.-H. Cheng, J.-F. Hwang, and P.~Yang, \emph{Existence and
uniqueness for
  {$p$}-area minimizers in the {H}eisenberg group}, Math. Ann. \textbf{337}
  (2007), no.~2, 253--293. \MR{MR2262784}

\bibitem{dgn}
D.~Danielli, N.~Garofalo, and D.-M. Nhieu, \emph{Sub-{R}iemannian
calculus on
  hypersurfaces in {C}arnot groups}, Adv. Math. \textbf{215} (2007), no.~1,
  292--378. \MR{MR2354992}

\bibitem{dgn3}
\bysame, \emph{A notable family of entire intrinsic minimal graphs
in the
  {H}eisenberg group which are not perimeter minimizing}, Amer. J. Math.
  \textbf{130} (2008), no.~2, 317--339. \MR{MR2405158}

\bibitem{dgnp-stable}
D.~Danielli, N.~Garofalo, D.~M. Nhieu, and S.~D. Pauls, \emph{{The
Bernstein
  problem for embedded surfaces in the Heisenberg group $\mathbb{H}^1$}},
  preprint. Revision of \emph{Stable complete embedded minimal surfaces in
  $\mathbb{H}^1$ with empty characteristic locus are vertical planes},
  arXiv:0903.4296.

\bibitem{dgnp}
\bysame, \emph{Instability of graphical strips and a positive
answer to the
  {B}ernstein problem in the {H}eisenberg group {$\Bbb H^1$}}, J. Differential
  Geom. \textbf{81} (2009), no.~2, 251--295. \MR{MR2472175}

\bibitem{d2}
M.~Derridj, \emph{Sur un th\'eor\`eme de traces}, Ann. Inst.
Fourier (Grenoble)
  \textbf{22} (1972), no.~2, 73--83. \MR{MR0343011 (49 \#7755)}

\bibitem{dcriem}
M.~P. do~Carmo, \emph{Riemannian geometry}, Mathematics: Theory \&
  Applications, Birkh\"auser Boston Inc., Boston, MA, 1992, Translated from the
  second Portuguese edition by Francis Flaherty. \MR{MR1138207 (92i:53001)}

\bibitem{dcp}
M.~P. do~Carmo and C.~K. Peng, \emph{Stable complete minimal
surfaces in {${\bf
  R}\sp{3}$} are planes}, Bull. Amer. Math. Soc. (N.S.) \textbf{1} (1979),
  no.~6, 903--906. \MR{MR546314 (80j:53012)}

\bibitem{evans}
L.~C. Evans and R.~F. Gariepy, \emph{Measure theory and fine
properties of
  functions}, Studies in Advanced Mathematics, CRC Press, Boca Raton, FL, 1992.
  \MR{MR1158660 (93f:28001)}

\bibitem{fc}
D.~Fischer-Colbrie, \emph{On complete minimal surfaces with finite
{M}orse
  index in three-manifolds}, Invent. Math. \textbf{82} (1985), no.~1, 121--132.
  \MR{MR808112 (87b:53090)}

\bibitem{fcs}
D.~Fischer-Colbrie and R.~Schoen, \emph{The structure of complete
stable
  minimal surfaces in {$3$}-manifolds of nonnegative scalar curvature}, Comm.
  Pure Appl. Math. \textbf{33} (1980), no.~2, 199--211. \MR{MR562550
  (81i:53044)}

\bibitem{fssc}
B.~Franchi, R.~Serapioni, and F.~Serra~Cassano,
\emph{Rectifiability and
  perimeter in the {H}eisenberg group}, Math. Ann. \textbf{321} (2001), no.~3,
  479--531. \MR{MR1871966 (2003g:49062)}

\bibitem{fsscadv}
\bysame, \emph{Regular submanifolds, graphs and area formula in
{H}eisenberg
  groups}, Adv. Math. \textbf{211} (2007), no.~1, 152--203. \MR{MR2313532
  (2008h:49030)}

\bibitem{gromov-cc}
M.~Gromov, \emph{Carnot-{C}arath\'eodory spaces seen from within},
  Sub-Riemannian geometry, Progress in Mathematics, vol. 144, Birkh\"auser,
  Basel, 1996, pp.~79--323. \MR{MR1421823 (2000f:53034)}

\bibitem{hp2}
R.~K. Hladky and S.~D. Pauls, \emph{{Variation of perimeter
measure in
  sub-Riemannian geometry}}, arXiv:math/0702237.

\bibitem{hp1}
\bysame, \emph{Constant mean curvature surfaces in
sub-{R}iemannian geometry},
  J. Differential Geom. \textbf{79} (2008), no.~1, 111--139. \MR{MR2401420}

\bibitem{hr2}
A.~Hurtado and C.~Rosales, \emph{Stable surfaces inside the
sub-{R}iemannian
  three-sphere}, in preparation.

\bibitem{hr1}
\bysame, \emph{Area-stationary surfaces inside the
sub-{R}iemannian
  three-sphere}, Math. Ann. \textbf{340} (2008), no.~3, 675--708. \MR{MR2358000
  (2008i:53038)}

\bibitem{msc}
R.~Monti and F.~Serra~Cassano, \emph{Surface measures in
  {C}arnot-{C}arath\'eodory spaces}, Calc. Var. Partial Differential Equations
  \textbf{13} (2001), no.~3, 339--376. \MR{MR1865002 (2002j:49052)}

\bibitem{mscv}
R.~Monti, F.~Serra~Cassano, and D.~Vittone, \emph{A negative
answer to the
  {B}ernstein problem for intrinsic graphs in the {H}eisenberg group},
  Bollettino dell unione matematica italiana (2008), no.~3, 709--728, ISSN
  1972-6724.

\bibitem{pauls-regularity}
Scott~D. Pauls, \emph{{$H$}-minimal graphs of low regularity in
{$\Bbb H\sp
  1$}}, Comment. Math. Helv. \textbf{81} (2006), no.~2, 337--381. \MR{MR2225631
  (2007g:53032)}

\bibitem{rr1}
M.~Ritor\'e and C.~Rosales, \emph{Rotationally invariant
hypersurfaces with
  constant mean curvature in the {H}eisenberg group {$\Bbb H\sp n$}}, J. Geom.
  Anal. \textbf{16} (2006), no.~4, 703--720. \MR{MR2271950}

\bibitem{rr2}
M.~Ritor{\'e} and C.~Rosales, \emph{Area-stationary surfaces in
the
  {H}eisenberg group {$\Bbb H\sp 1$}}, Adv. Math. \textbf{219} (2008), no.~2,
  633--671. \MR{MR2435652}

\bibitem{r2}
Manuel Ritor{\'e}, \emph{Examples of area-minimizing surfaces in
the
  sub-{R}iemannian {H}eisenberg group {$\Bbb H1$} with low regularity}, Calc.
  Var. Partial Differential Equations \textbf{34} (2009), no.~2, 179--192.
  \MR{MR2448649 (2009h:53062)}

\bibitem{ros}
A.~Ros, \emph{One-sided complete stable minimal surfaces}, J.
Differential
  Geom. \textbf{74} (2006), no.~1, 69--92. \MR{MR2260928 (2007g:53008)}

\bibitem{rcmc}
C.~Rosales, \emph{Complete stable surfaces under a volume
constraint in the
  first {H}eisenberg group}, in preparation.

\bibitem{simon}
L.~Simon, \emph{Lectures on geometric measure theory}, Proceedings
of the
  Centre for Mathematical Analysis, Australian National University, vol.~3,
  Australian National University Centre for Mathematical Analysis, Canberra,
  1983. \MR{MR756417 (87a:49001)}

\end{thebibliography}
\end{document}